\documentclass[12pt]{article}
\usepackage{amssymb,amsmath, txfonts}
\begin{document}

\def\abstractname{\bf Abstract}
\def\dfrac{\displaystyle\frac}
\let\oldsection\section
\renewcommand\section{\setcounter{equation}{0}\oldsection}
\renewcommand\thesection{\arabic{section}}
\renewcommand\theequation{\thesection.\arabic{equation}}
\newtheorem{theorem}{\indent Theorem}[section]
\newtheorem{lemma}{\indent Lemma}[section]
\newtheorem{proposition}{\indent Proposition}[section]
\newtheorem{definition}{\indent Definition}[section]
\newtheorem{remark}{\indent Remark}[section]
\newtheorem{corollary}{\indent Corollary}[section]
\def\pd#1#2{\displaystyle\frac{\partial#1}{\partial#2}}
\def\d#1{\displaystyle\frac{d#1}{dt}}

\title{\LARGE\bf Global bounded solution  in three-dimensional chemotaxis-Stokes model with arbitrary porous medium slow diffusion
\thanks{This work is supported by NSFC(11871230), Guangdong Basic and Applied Basic Research Foundation(2020B1515310013).}
\\
\author{Chunhua Jin$^a$\thanks{
Corresponding author. Email: {\tt jinchhua@126.com}}
\\
\small \it{$^a$School of Mathematical Sciences, South China Normal University, }
\\
\small \it{Guangzhou, 510631, China}
}}

\date{}

\maketitle

\begin{abstract}
In this paper, we study the following chemotaxis-Stokes  model with porous medium slow diffusion
\begin{align*}
\left\{
\begin{aligned}
&n_t+u\cdot\nabla n=\Delta n^m-\chi\nabla\cdot(n\nabla c),
\\
&c_t+u\cdot\nabla c-\Delta c=-cn,
\\
&u_t+\nabla \pi=\Delta u+n\nabla\varphi,
\\
& {\rm div} u=0,
\\
\end{aligned}\right.
\end{align*}
in a bounded domain $\Omega\subset \mathbb R^3$ with zero-flux boundary conditions and no-slip boundary condition.
In recent ten years, many efforts have been made to find the global bounded solutions of chemotaxis-Stokes system in three dimensional space.
Although some important progress  has been carried out in some papers \cite{DLM, TW2, W2, W3}, as mentioned by some authors, the question of identifying an optimal condition on $m\ge 1$ ensuring global boundedness in the three-dimensional framework remains an open challenge.
In the present paper,
we  put forward  a new estimation technique, completely proved the existence of global bounded solutions for arbitrary slow diffusion case $m>1$,
and partially answered the open problem proposed by Winkler.
\end{abstract}


%%  92C17  35K40  35K65

\section{Introduction}
In this paper, we consider the following typical chemotaxis-Stokes system
\allowdisplaybreaks
\begin{align}
\label{1-1}\left\{
\begin{aligned}
&n_t+u\cdot\nabla n=\Delta n^m-\chi\nabla\cdot(n\nabla c), &&    (x,t) \in Q,
\\
&c_t+u\cdot\nabla c-\Delta c=-cn, &&  (x,t) \in Q,
\\
&u_t+\nabla \pi=\Delta u+n\nabla\varphi, &&    (x,t) \in Q,
\\
& {\rm div} u=0, && (x,t)\in Q,
\\
&\left.(\nabla n^m-\chi n\nabla c)\cdot{\nu}\right|_{\partial\Omega}
=\left.\frac{\partial c}{\partial {\nu}}\right|_{\partial\Omega}=0,   \quad u|_{\partial\Omega}=0,
\\
& n(x,0)=n_0(x),  c(x,0)=v_0(x), u(x,0)=u_0(x), &&  x\in\Omega,
\end{aligned}\right.
\end{align}
where $m>1$, $Q=\Omega\times \mathbb R^+$,
$\Omega\subset \mathbb R^3$ is a bounded domain, and the boundary $\partial\Omega$ is appropriately smooth,
$n$, $c$ represent the bacteria cell density, the oxygen concentration respectively, $u$, $\pi$ are
the fluid velocity and the associated pressure, $\chi>0$ is the sensitivity coefficient of aggregation induced by the  concentration changes of  oxygen,
$-cn$ is the consumption term of oxygen,  that is  more bacteria,  more oxygen are consumed, the fluid couples to $n$
and $c$ through transports $u\cdot\nabla n$, $u\cdot\nabla c$ and the gravitational potential
$n\nabla\varphi$.

The chemotaxis-fluid model was initially introduced by Tuval, Goldstein et. al  \cite{TC}  in 2005,
which is a coupled system of  chemotaxis-Navier-Stokes equations as follows
\begin{align}\label{1-3}\left\{
\begin{aligned}
&n_t+u\cdot\nabla n=\Delta n^m-\chi\nabla\cdot(n\nabla c),
\\
&c_t+u\cdot\nabla c-\Delta c=-cn,
\\
&u_t+\tau u\cdot\nabla u+\nabla \pi=\Delta u+n\nabla\varphi,
\\
& {\rm div} u=0.
\end{aligned}\right.
\end{align}
This model is derived from the results of an experimental observation, which describes the dynamics of bacterial swimming and oxygen transport near
contact lines. In this model, the fluid motion is governed by the full Navier-Stokes equations ($\tau=1$).
However, it is well known that
the existence of smooth solutions for the 3-D incompressible Navier-Stokes equations has always been an open problem.
So, some mathematicians turn their attention to a simpler model, that is the case $\tau=0$.
In fact, when the fluid
motion is slow, the inertial force is much less than the viscous force,  for this case, the Navier-Stokes equations indeed can be replaced by Stokes equations
by ignoring the convective term \cite{KP}.

The chemotaxis-(Navier)-Stokes model has attracted many researcher's attention since it was proposed.
For the two dimensional case of \eqref{1-3}, the global solvability and boundedness of weak solutions are established completely
for any $m>1$ in \cite{TW1}.
While in three dimensional space, the research of \eqref{1-3} ($\tau=0$) is rather tortuous. The first effort to this 3-D problem is due to the
work of Di Francesco et al. \cite{DLM}, in which, they established the existence of global bounded
weak solutions for $m$ in some finite interval, namely  $m\in \left(\frac{7+\sqrt{217}}{12}, 2\right]$ (approximating to  $(1.81, 2]$);
Then, Tao and Winkler \cite{TW2}, in 2013, who obtained the  global existence of locally
bounded weak solutions with $m$ belonging to the infinite interval $(\frac87, +\infty)$.
Afterwards, Winkler \cite{W2} supplemented the uniform boundedness of solutions
for the case $m>\frac 76$; Recently, Winkler \cite{W3} further improved this result to the case $m>\frac 98$.
However, as mentioned by Winkler, the question of identifying an optimal condition on
$m\ge 1$ ensuring global boundedness in the three-dimensional version of chemotaxis-Stokes system  remains an open challenge.

In the present paper, we  pay our attention to the global existence  and uniform boundedness  of weak solutions
for the system \eqref{1-1}. The new ideas lies in that
firstly, we make full use of the consumption term, even if when do gradient dependent estimation on $c$,
we bring out some good terms from them.
Secondly, we proposed a new estimation technique on $\frac{|\nabla n|^2}{n}$ using some very meticulous analysis, see Lemma \ref{lem3-5}.
Finally, we show that for any $m>1$,
this problem admits a global  bounded weak solution
for any large initial datum.  Our research for this system improved the result of \cite{DLM, TW2, W2, W3}.
Moreover, we also improve the regularity for $m\le \frac54$, and proved that the obtained solutions are strong solutions.

Before going further, we first give the assumption on initial data.

Assumptions:
\begin{align*}
\qquad\qquad\qquad\qquad
\left\{
\begin{aligned}
& n_0\in L^\infty(\Omega), \nabla \sqrt n_0\in L^{2}(\Omega), \sqrt c_0\in L^{2}(\Omega), c_0\in W^{2,p}(\Omega)\  \text{for any $p>1$},
\\
&u_0\in  W^{2,p}(\Omega) \  \text{for any $p>1$},  \varphi\in W^{1,\infty}(\Omega),
\\
& n_0, c_0\ge 0.
\end{aligned}\right.\quad (H)
\end{align*}
From (H), by Sobolev imbedding theorem, it is easy to see that $c_0, u_0\in  W^{1,\infty}(\Omega)$ and $A^{\beta}u_0\in L^p(\Omega)$ for any $\beta\in(\frac12, 1)$, $p>1$.

The main results:
\begin{theorem}
\label{thm-1}
Assume $(H)$,  $m>1$.
Then the problem \eqref{1-1} admits a  global bounded weak solution $(n,c,u,\pi)$
with  $n\in\mathcal X_1, c\in \mathcal X_2$, $u\in \mathcal X_3$, $\pi\in  \mathcal X_4$,
where
$$
\mathcal X_1=\{n\in L^\infty(\Omega\times \mathbb R^+); \nabla n^m\in L^\infty(\mathbb R^+;L^2(\Omega)),
\left(n^{\frac{m+1}2}\right)_t, \nabla n^{\frac{m}2} \in L^2_{loc}([0,\infty); L^2(\Omega))\},
$$
$$
\mathcal X_2=\{c\in L^\infty(\mathbb R^+; W^{1,\infty}(\Omega)); c_t,  D^2c\in  L^p_{loc}([0,\infty);L^p(\Omega))
\  \text{for any }\  p>1\},
$$
$$
\mathcal X_3=\left\{u\in L^\infty(\mathbb R^+; W^{1,\infty}(\Omega)); u_t,  D^2u\in  L^p_{loc}([0,\infty);L^p(\Omega))
\  \text{for any }\  p>1 \right\},
$$
$$
\mathcal X_4=\{\pi; \nabla\pi\in L^p_{loc}([0,\infty);L^p(\Omega))
\  \text{for any }\  p>1\},
$$
such that
\begin{align}
\label{1-3}
&\sup_{t\in(0,+\infty)}\left(\|n(\cdot,t)\|_{L^\infty}+\|c(\cdot,t)\|_{W^{1,\infty}}+\|u(\cdot,t)\|_{W^{1,\infty}}
+\|\nabla n^m(\cdot,t)\|_{L^2}\right)\le M_1,
\\
\label{1-4}
&\sup_{t\in (0,+\infty)}\left(\|(n^{\frac{m+1}2})_t\|_{L^2(Q_1(t))}+\|\nabla n^{\frac{m}2}\|_{L^2(Q_1(t))}+
\|c\|_{W_p^{2,1}(Q_1(t))}+\|u\|_{W_p^{2,1}(Q_1(t))}+\|\nabla \pi\|_{L^p(Q_1(t))}
\right)\le M_2,
\end{align}
for any $p>1$. Here $Q_1(t)=\Omega\times (t, t+1)$,
$M_i$ $(i=1,2)$  depend  only on $\chi$, $\Omega$, $m$,  $p$, $n_0$, $c_0$, $u_0$.
\end{theorem}
In particular, when $m\le \frac54$, the solution we obtained is strong solution.
That is
\begin{theorem}
\label{thm-2}
Assume $(H)$,  $1<m\le \frac54$.
Then the problem \eqref{1-1} admits a  global bounded strong  solution $(n,c,u,\pi)$
with  $n\in\mathcal D_1, c\in \mathcal X_2$, $u\in \mathcal X_3$, $\pi\in  \mathcal X_4$, such that
\begin{align*}
&\iint_{Q_T}\left(\frac{\partial n}{\partial t}+u\cdot\nabla n-\Delta n^m+\chi\nabla\cdot (n\nabla c)\right)\phi dxdt=0, \qquad \text{for any $\phi\in L^2(Q_T)$,}
\\
&\iint_{Q_T}\left(\frac{\partial c}{\partial t}+u\cdot\nabla c-\Delta c+n c\right)\psi dxdt=0,
\qquad \text{for any $\psi\in L^q(Q_T)$,}
\\
&\iint_{Q_T}\left(\frac{\partial u}{\partial t}-\Delta u+\nabla \pi-n\nabla \varphi\right)\Phi dxdt=0,
\qquad \text{for any $\Phi\in L^q(Q_T)$,}
\end{align*}
where
$$
\mathcal D_1=\{n\in L^\infty(\Omega\times \mathbb R^+); \nabla\sqrt n\in L^\infty(\mathbb R^+;L^2(\Omega)),
n_t, \Delta n^m, n^{\frac m2-2}|\nabla n|^2, n^{\frac12}| D^2 n^{\frac{m-1}2}| \in L^2_{loc}([0,\infty); L^2(\Omega))\},
$$
and
\begin{align}
\label{1-5}
&\sup_{t\in(0,+\infty)}\left(\|n(\cdot,t)\|_{L^\infty}+\|c(\cdot,t)\|_{W^{1,\infty}}+\|u(\cdot,t)\|_{W^{1,\infty}}
+\|\nabla {\sqrt n}(\cdot,t)\|_{L^2}\right)\le M_1,
\\
&\sup_{t\in (0,+\infty)}\left(\|n_t\|_{L^2(Q_1(t))}+\|\Delta n^m\|_{L^2(Q_1(t))}+\|n^{\frac{m}4-1}\nabla n\|_{L^4(Q_1(t))}
+\|n^{\frac12} D^2 n^{\frac{m-1}2}\|_{L^2(Q_1(t))}\right.\nonumber
\\
\label{1-6}
&\qquad \left.+\|c\|_{W_p^{2,1}(Q_1(t))}+\|u\|_{W_p^{2,1}(Q_1(t))}+\|\nabla \pi\|_{L^p(Q_1(t))}
\right)\le M_2, \quad \forall  \  p>1.
\end{align}
Here $Q_1(t)=\Omega\times (t, t+1)$,
$M_i$ $(i=1,2)$ are constants depending only on $\chi$, $\Omega$, $m$,  $p$, $n_0$, $c_0$, $u_0$.
\end{theorem}

\section{Preliminaries}

We first give some notations, which will be used throughout this paper.

\medskip

{\bf Notations:}
$\|\cdot\|_{L^p}=\|\cdot\|_{L^p(\Omega)}$,
$Q_{1}(t)=\Omega\times(t, t+1)$,  $Q_T:=Q_T(0)=\Omega\times (0, T)$.

\medskip

Before going further, we list some important lemmas, which will be used throughout this paper.
By \cite{BB} or \cite{LM}, we first give the following lemma.
\begin{lemma}
\label{lem2-1}
Given an $f\in W^{k,p}(\Omega)$ ($k\ge 0$ an integer) and a $g\in W^{k+1-\frac1p, p}(\partial\Omega)$
such that
$$
\int_\Omega fdx=\int_{\partial\Omega} gdS,
$$
there exists a $u\in W^{k+2, p}(\Omega)$ satisfying
\begin{align*}\left\{
\begin{aligned}
&\Delta u=f, \qquad x\in \Omega
\\
&\left.\frac{\partial u}{\partial{\bf n}}\right|_{\partial\Omega}=g,
\end{aligned}\right.
\end{align*}
In addition
$$
\|\nabla u\|_{W^{k+1,p}(\Omega)}\le C\left(\|f\|_{W^{k,p}(\Omega)}+\|g\|_{W^{k+1-\frac1p, p}(\partial\Omega)}\right).
$$
\end{lemma}

\medskip

From this lemma, it is easy to get that
$$
\|\nabla u\|_{W^{k+1,p}(\Omega)}\le C\|\Delta u\|_{W^{k,p}(\Omega)}.
$$
if $\left.\frac{\partial u}{\partial{\bf n}}\right|_{\partial\Omega}=0$ by taking $f=\Delta u$, $g=0$ in Lemma \ref{lem2-1}.

\medskip

By \cite{J1}, we have the following  $L^p$-theory of linear parabolic equations.

\begin{lemma}
\label{lem2-3}
Assume that $u_0\in W^{2,p}(\Omega)$, and  $f\in L_{loc}^p((0,+\infty);L^p(\Omega))$
with
$$
\sup_{t\in(\tau, +\infty)}\int_{t-\tau}^{t}\|f\|_{L^p}^pds\le A,
$$
where $\tau>0$ is a fixed constant.
Then the following problem
\begin{align}
\label{LP}\left\{
\begin{aligned}
&u_t-\Delta u+u=f(x,t),
\\
&\left.\frac{\partial u}{\partial{\bf n}}\right|_{\partial\Omega}=0,
\\
&u(x,0)=u_0(x)
\end{aligned}\right.
\end{align}
admits a unique solution   $u\in  L_{loc}^p((0,+\infty);W^{2,p}(\Omega))$, $u_t\in L_{loc}^p((0,+\infty);L^p(\Omega))$
with
\begin{equation}
\label{LPE}
\sup_{t\in(\tau, +\infty)}\int_{t-\tau}^{t}(\|u\|_{W^{2,p}}^p+\|u_t\|_{L^p}^p)ds\le
AM\frac{e^{p\tau}}{e^{\frac p2\tau}-1}+Me^{\frac p2\tau}\|u_0\|_{W^{2,p}}^p,
\end{equation}
where $M$ is a constant independent of $\tau$.
\end{lemma}

By \cite{DG, W1, MS}, we have the following three lemmas.
\begin{lemma}
\label{lem2-4}
Suppose that $h\in C^2(\mathbb R)$, then for all $\varphi\in C^2(\overline \Omega)$
fulfilling $\frac{\partial\varphi}{\partial {\bf n}}=0$ on $\partial\Omega$, we have
\begin{align*}
&\int_\Omega h'(\varphi)|\nabla\varphi|^2\Delta\varphi dx+\frac23\int_\Omega
h(\varphi)|\Delta\varphi|^2 dx\nonumber
\\
=&\frac23\int_\Omega h(\varphi)| D^2\varphi|^2 dx
-\frac13\int_\Omega h''(\varphi)|\nabla\varphi|^4dx-\frac13\int_{\partial\Omega}
h(\varphi)\frac{\partial|\nabla\varphi|^2}{\partial {\bf n}}ds,
\end{align*}
where $|D^2\varphi|^2=\sum_{i,j=1}^n |D_{ij}\varphi|^2$.
\end{lemma}
\begin{lemma}
\label{lem2-5}
Suppose that $h\in C^1(\mathbb R^+)$ is positive, and let $\theta(s)=\int_1^s\frac{1}{h(t)}dt$ for $s>0$,
then for all $\varphi\in C^2(\overline \Omega)$
fulfilling $\frac{\partial\varphi}{\partial {\bf n}}\Big|_{\partial\Omega}=0$, the following inequality holds
$$
\int_\Omega\frac{h'(\varphi)}{h^3(\varphi)}|\nabla \varphi|^4 dx\le (2+\sqrt N)^2\int_\Omega \frac{h(\varphi)}{h'(\varphi)}| D^2\theta(\varphi)|^2dx.
$$
\end{lemma}

\begin{lemma}
\label{lem2-6}
Assume that $\Omega$ is bounded and let $\omega\in C^2(\overline\Omega)$ satisfy
$\frac{\partial\omega}{\partial\nu}\Big|_{\partial\Omega}=0$. Then we have
$$
\frac{\partial|\nabla\omega|^2}{\partial\nu}\le 2\kappa|\nabla\omega|^2 \quad \text{on} \ \partial\Omega,
$$
where  $\kappa>0$ is an upper bound for the curvatures of $\Omega$.
\end{lemma}

\section{Global Classical Solution and Uniform Energy Estimates to the Regularized Problem}
To study the existence of solutions to the system \eqref{1-1}, we
consider the following approximate problems
\begin{align}
\label{3-1}\left\{
\begin{aligned}
&n_{\varepsilon t}+u_{\varepsilon}\cdot\nabla n_{\varepsilon}=m\nabla\cdot((\varepsilon +n_{\varepsilon})^{m-1}\nabla n_{\varepsilon} )-\chi
\nabla\cdot\left(n_{\varepsilon}\nabla c_{\varepsilon}\right),  &&  (x,t) \in Q,
\\
&c_{\varepsilon t}+u_{\varepsilon}\cdot\nabla c_{\varepsilon}-\Delta c_{\varepsilon}=-c_{\varepsilon}n_{\varepsilon}, &&  (x,t) \in Q,
\\
&u_{\varepsilon t}+\nabla \pi_{\varepsilon}=\Delta u_{\varepsilon}+n_{\varepsilon}\nabla\varphi_{\varepsilon}, &&    (x,t) \in Q,
\\
& {\rm div} u_{\varepsilon}=0, &&   (x,t) \in Q,
\\
&\left.\frac{\partial n_{\varepsilon}}{\partial {\bf n}}\right|_{\partial\Omega}
=\left.\frac{\partial c_{\varepsilon}}{\partial {\bf n}}\right|_{\partial\Omega}=0,
\quad u_{\varepsilon}|_{\partial\Omega}=0, &&
\\
& n_{\varepsilon}(x,0)=n_{\varepsilon 0}(x)\ge 0,  c_{\varepsilon}(x,0)=c_{\varepsilon 0}(x)\ge 0,
u_{\varepsilon}(x,0)=u_{\varepsilon 0}(x), &&  x\in\Omega
\end{aligned}\right.
\end{align}
for any $\varepsilon\in(0,1)$,
where $\varphi_\varepsilon\in C^{1+\alpha,\frac\alpha 2}(\overline{\Omega}\times[0,+\infty))$,
$n_{\varepsilon 0}, c_{\varepsilon 0}, u_{\varepsilon 0} \in C^{2+\alpha}(\overline\Omega)$ with
$\left.\frac{\partial n_{\varepsilon 0}}{\partial{\bf n}}\right|_{\partial\Omega}=
0, \left.\frac{\partial c_{\varepsilon 0}}{\partial{\bf n}}\right|_{\partial\Omega}=0$,
$\left.u_{\varepsilon 0}\right|_{\partial\Omega}=0$,
$$
n_{\varepsilon0} \rightarrow n_0, c_{\varepsilon0}\rightarrow c_0,
u_{\varepsilon0}\rightarrow u_0, \nabla\varphi_\varepsilon\rightarrow\nabla\varphi,\quad
\text{strongly in $L^p$ for any $p>1$},
$$
and
\begin{align*}
&\|n_{\varepsilon0}\|_{L^\infty}+\left\|\nabla \sqrt n_{\varepsilon0}\right\|_{L^2}+\|\nabla \sqrt c_{\varepsilon0}\|_{L^2}
+\|c_{\varepsilon0}\|_{W^{1,\infty}}+\| D^2 c_{\varepsilon0}\|_{L^p}+\|u_{\varepsilon0}\|_{W^{2,p}}
+\|u_{\varepsilon0}\|_{L^\infty}+\|\nabla\varphi_\varepsilon\|_{L^\infty}
\\
\le &2\left(\|n_{0}\|_{L^\infty}+\left\|\nabla \sqrt n_0\right\|_{L^2}++\left\|\nabla \sqrt c_0\right\|_{L^2}
+\|c_{0}\|_{W^{1,\infty}}+\| D^2 c_0\|_{L^p}+\|u_{0}\|_{W^{2,p}}
+\|u_{0}\|_{L^\infty}+\|\nabla\varphi\|_{L^\infty}\right).
\end{align*}
Let $A$  be the Stokes operator, that is $Aw:=-P\Delta w$,
$P:L^r(\Omega)\rightarrow L^r_\sigma(\Omega)$ is the Helmholtz projection  \cite{G1}.
$A$ generates a bounded
analytic semigroup $\{e^{-tA}\}_{t\ge 0}$ on $L_\sigma^r$, and the solution $u$ of \eqref{1-1}
can be expressed as
\begin{equation}
\label{3-2}
u=e^{-tA}u_0+\int_{0}^t e^{-(t-s)A}P(n(s)\nabla\varphi(s))ds.
\end{equation}
For more details of Stokes operator, we refer to \cite{KY}. The following local existence result of classical solution
is well known,
We first state the local existence result of classical solution to \eqref{3-1} as follows, which is well known,
see for example \cite{TW2, J1}.

\begin{lemma}
\label{lem3-1}
Assume that $m>1$,   there exists $T_{\max}\in(0, +\infty]$ such that
the problem \eqref{3-1} admits a unique classical solution
$(n_\varepsilon, c_\varepsilon, u_\varepsilon, \pi_\varepsilon)\in C^{2+\alpha, 1+\alpha/2}(\overline\Omega \times(0, T_{\max}))$ with
$$
n_\varepsilon\ge 0, \quad c_\varepsilon\ge 0,  \ \text{for all}\  (x, t)\in \Omega\times (0, T_{\max}),
$$
such that
either $T_{max} = \infty$, or
$$
\limsup_{t\nearrow T_{\max}}\left(\|n_\varepsilon(\cdot, t)\|_{L^\infty}+\|c_\varepsilon(\cdot, t)\|_{W^{1,\infty}}\right)
= \infty.
$$
\end{lemma}

In the following lemmas of this section,  we denote $\tau=\min\{1, \frac{T_{\max}}2\}$.
Noticing that $\tau\le 1$, if $\tau<1$, it implies that $T_{\max}<2$, then the integral $\int_{t-\tau}^t\cdots$ of these lemmas in this section can be replaced with $\int_0^{T_{\max}}\cdots$. For simplicity, throughout this paper, we let $C$, $C_i$, $\tilde C$, $\hat C$  et.al denote some different constants, which are independent of $\varepsilon$ and $T_{\max}$, and $\tau$, and if no special explanations,
which depend at most on $n_0$, $c_0$, $u_0$, $m$, $\chi$, $\|\nabla\varphi\|_{L^\infty}$ and $\Omega$.

By maximum principle, it is easy to see that $c_\varepsilon$ is bounded,
and by  a direct integration for the first equation of \eqref{3-1}, it is also easy to
obtain the $L^1$-norm of $n_\varepsilon$, that is,
\begin{lemma}
\label{lem3-2}
Let $(n_\varepsilon, c_\varepsilon, u_\varepsilon, \pi_\varepsilon)$ be the classical
solution of \eqref{3-1} in $(0, T_{\max})$, then we have
\begin{align}
\label{3-3}
&\sup_{t\in(0,T_{\max})}\|c_\varepsilon(\cdot,t)\|_{L^\infty}\le C_1,
\\
\label{3-4}
&\sup_{t\in(0,T_{\max})}\int_\Omega n_\varepsilon dx\le C_2,
\end{align}
where $C_1$, $C_2$ are independent of $\varepsilon$.
\end{lemma}

Next, we show the following lemma.
\begin{lemma}
\label{lem3-3}
Let $(n_\varepsilon, c_\varepsilon, u_\varepsilon, \pi_\varepsilon)$ be the classical
solution of \eqref{3-1}. Then
\begin{align}
&\sup_{t\in(0,T_{\max})}\int_\Omega\left(\frac{|\nabla c_\varepsilon|^2}{c_\varepsilon}+
n_{\varepsilon}\ln n_{\varepsilon}+ |u_\varepsilon|^2\right)dx+\sup_{t\in(\tau,T_{\max})}\int_{t-\tau}^{t}
\int_\Omega \left(c_{\varepsilon}| D^2\ln c_{\varepsilon}|^2+\frac{| D^2 c_{\varepsilon}|^2}{c_{\varepsilon}}+
\frac{|\nabla c_\varepsilon|^4}{|c_\varepsilon|^3}\right)dx
\nonumber
\\
\label{3-5}
&+\sup_{t\in(\tau,T_{\max})}\int_{t-\tau}^{t}
\int_\Omega \left(\frac{n_\varepsilon}{c_\varepsilon}|\nabla c_\varepsilon|^2 +\frac{n_{\varepsilon}+\varepsilon}{n_{\varepsilon}}
|\nabla (n_{\varepsilon}+\varepsilon)^{\frac m2}|^2+ |\nabla u_\varepsilon|^2\right)dx\le C,
\end{align}
where $C$ is independent of $\varepsilon$.
\end{lemma}

{\it\bfseries Proof.} Using the second equation of \eqref{3-1}, we see that
\begin{align}
&\frac12\frac{d}{dt}\int_\Omega\frac{|\nabla c_\varepsilon|^2}{c_\varepsilon}dx
=\int_\Omega\frac{\nabla c_\varepsilon}{c_\varepsilon}\nabla c_{\varepsilon t}dx
-\frac12\int_\Omega\frac{|\nabla c_\varepsilon|^2}{c_\varepsilon^2}c_{\varepsilon t}dx
\nonumber
\\
&=-\int_\Omega c_{\varepsilon t}\left(\frac{\Delta c_{\varepsilon}}{c_{\varepsilon}}
-\frac{|\nabla c_\varepsilon|^2}{c_\varepsilon^2}\right)
-\frac12\int_\Omega\frac{|\nabla c_\varepsilon|^2}{c_\varepsilon^2}c_{\varepsilon t}dx
\nonumber
\\
&=-\int_\Omega \frac{\Delta c_{\varepsilon}}{c_{\varepsilon}} \left(
\Delta c_{\varepsilon}- c_{\varepsilon}n_{\varepsilon}-u_\varepsilon\cdot\nabla c_\varepsilon\right)dx
+\frac12\int_\Omega \frac{|\nabla c_\varepsilon|^2}{c_\varepsilon^2}
\left(\Delta c_{\varepsilon}- c_{\varepsilon}n_{\varepsilon}-u_\varepsilon\cdot\nabla c_\varepsilon)\right)dx
\nonumber
\\
&=-\int_\Omega \frac{|\Delta c_{\varepsilon}|^2}{c_{\varepsilon}}
+\frac12\int_\Omega \frac{|\nabla c_\varepsilon|^2}{c_\varepsilon^2}\Delta c_{\varepsilon} dx
-\int_\Omega \left(\nabla n_{\varepsilon}\nabla c_{\varepsilon}
+\frac{1}2\frac{n_\varepsilon}{c_\varepsilon}|\nabla c_\varepsilon|^2\right)dx\nonumber
\\
\label{3-6}
&+\int_\Omega\frac{u_\varepsilon\cdot\nabla c_\varepsilon}{c_\varepsilon}\Delta c_\varepsilon dx
-\frac12\int_\Omega u_\varepsilon\cdot\nabla c_\varepsilon\frac{|\nabla c_\varepsilon|^2}{c_\varepsilon^2}dx.
\end{align}
By a direct calculation, we get that
\begin{align}
\int_\Omega c_{\varepsilon}| D^2\ln c_{\varepsilon}|^2 dx &=\int_\Omega
\left(\frac{| D^2c_{\varepsilon}|^2}{c_{\varepsilon}}+\frac{|\nabla c_\varepsilon|^4}
{c_\varepsilon^3}-2\frac{\nabla c_\varepsilon \cdot  D^2
c_\varepsilon\cdot\nabla c_\varepsilon}{c_\varepsilon^2}\right) dx\nonumber
\\
&=\int_\Omega
\left(\frac{| D^2c_{\varepsilon}|^2}{c_{\varepsilon}}+\frac{|\nabla c_\varepsilon|^4}
{c_\varepsilon^3}-\frac{\nabla c_\varepsilon\nabla |\nabla c_\varepsilon|^2}{c_\varepsilon^2}
\right)dx\nonumber
\\
&=\int_\Omega
\left(\frac{| D^2c_{\varepsilon}|^2}{c_{\varepsilon}}+\frac{|\nabla c_\varepsilon|^4}
{c_\varepsilon^3}+\frac{|\nabla c_\varepsilon|^2}{c_\varepsilon^2}\Delta c_\varepsilon
-2\frac{|\nabla c_\varepsilon|^4}{c_\varepsilon^3}\right)dx\nonumber
\\
\label{add3-7}
&=\int_\Omega
\left(\frac{| D^2c_{\varepsilon}|^2}{c_{\varepsilon}}+\frac{|\nabla c_\varepsilon|^2}{c_\varepsilon^2}\Delta c_\varepsilon
-\frac{|\nabla c_\varepsilon|^4}{c_\varepsilon^3}\right)dx.
\end{align}
Recalling strong maximum principle, $c_\varepsilon(x,t)>0$ for any $x\in\Omega$, $t>0$.
Using Lemma \ref{lem2-4} and combining with \eqref{add3-7},  we derive that
\begin{align*}
&-\int_\Omega \frac{|\Delta c_{\varepsilon}|^2}{c_{\varepsilon}}
+\frac12\int_\Omega \frac{|\nabla c_\varepsilon|^2}{c_\varepsilon^2}\Delta c_{\varepsilon} dx
\\
=&-\int_\Omega\frac{|\nabla c_\varepsilon|^2}{c_\varepsilon^2}\Delta c_{\varepsilon} dx
-\int_\Omega\frac1{c_\varepsilon}| D^2c_\varepsilon|^2 +\frac{1}{c_\varepsilon^3}
|\nabla c_\varepsilon|^4 dx+\frac12\int_{\partial\Omega}\frac{1}{c_\varepsilon}
\frac{\partial}{\partial{\bf n}}|\nabla c_\varepsilon|^2 ds
\\
=&-\int_\Omega c_{\varepsilon}| D^2\ln c_{\varepsilon}|^2 dx
+\frac12\int_{\partial\Omega}\frac{1}{c_\varepsilon}
\frac{\partial}{\partial{\bf n}}|\nabla c_\varepsilon|^2 ds.
\end{align*}
Substituting this equality into \eqref{3-6} yields
\begin{align}
&\frac12\frac{d}{dt}\int_\Omega\frac{|\nabla c_\varepsilon|^2}{c_\varepsilon}dx
+\int_\Omega c_{\varepsilon}| D^2\ln c_{\varepsilon}|^2 dx+
\frac{1}2\int_\Omega
\frac{n_\varepsilon}{c_\varepsilon}|\nabla c_\varepsilon|^2dx
+\int_\Omega
\nabla n_{\varepsilon}\nabla c_{\varepsilon}dx \nonumber
\\
&=\frac12\int_{\partial\Omega}\frac{1}{c_\varepsilon}
\frac{\partial}{\partial{\bf n}}|\nabla c_\varepsilon|^2 ds
+\int_\Omega\frac{u_\varepsilon\cdot\nabla c_\varepsilon}{c_\varepsilon}\Delta c_\varepsilon dx
-\frac12\int_\Omega u_\varepsilon\cdot\nabla c_\varepsilon\frac{|\nabla c_\varepsilon|^2}{c_\varepsilon^2}dx
\nonumber
\\
\label{3-7}
&=\frac12\int_{\partial\Omega}\frac{1}{c_\varepsilon}
\frac{\partial}{\partial{\bf n}}|\nabla c_\varepsilon|^2 ds
-\int_\Omega\frac{\nabla c_\varepsilon\cdot\nabla u_\varepsilon\cdot\nabla c_\varepsilon}{c_\varepsilon}dx.
\end{align}
Noticing that ${\rm div}u_\varepsilon=0$, then for $n_\varepsilon$, multiplying the first equation of \eqref{3-1} by $1+\ln n_{\varepsilon}$, we get
that
\begin{align}\label{3-8}
\frac{d}{dt}\int_\Omega n_{\varepsilon}\ln n_{\varepsilon}dx
+\frac{4}{m}\int_\Omega\frac{n_{\varepsilon}+\varepsilon}{n_{\varepsilon}}
|\nabla (n_{\varepsilon}+\varepsilon)^{\frac m2}|^2dx
=\chi\int_\Omega \nabla n_{\varepsilon}\nabla c_{\varepsilon} dx.
\end{align}
Combining \eqref{3-7}, and \eqref{3-8}, and using \eqref{3-3}, we obtain that
\begin{align}
&\frac{d}{dt}\int_\Omega\left(\frac12\frac{|\nabla c_\varepsilon|^2}{c_\varepsilon}+
\frac{1}{\chi}n_{\varepsilon}\ln n_{\varepsilon}\right)dx
+\int_\Omega \left(c_{\varepsilon}| D^2\ln c_{\varepsilon}|^2+\frac12\frac{n_\varepsilon}{c_\varepsilon}|\nabla c_\varepsilon|^2 +\frac{4}{m\chi}\frac{n_{\varepsilon}+\varepsilon}{n_{\varepsilon}}
|\nabla (n_{\varepsilon}+\varepsilon)^{\frac m2}|^2\right)dx\nonumber
\\
&=\frac12\int_{\partial\Omega}\frac{1}{c_\varepsilon}
\frac{\partial}{\partial{\bf n}}|\nabla c_\varepsilon|^2 ds
-\int_\Omega\frac{\nabla c_\varepsilon\cdot\nabla u_\varepsilon\cdot\nabla c_\varepsilon}{c_\varepsilon}dx\nonumber
\\
\label{3-9}
 &\le \frac12\int_{\partial\Omega}\frac{1}{c_\varepsilon}
\frac{\partial}{\partial{\bf n}}|\nabla c_\varepsilon|^2 ds+
\eta\int_\Omega \frac{|\nabla c_\varepsilon|^4}{c_\varepsilon^3}dx+C_\eta\int_\Omega
|\nabla u_\varepsilon|^2 dx
\end{align}
for any sufficiently small $\eta>0$, where $C_\eta$ depends on $\eta$.
By the boundary trace embedding theorem \cite{A}  and Lemma \ref{lem2-6},  we see that
for any sufficiently small $\eta_1>0$,
\begin{align}
&\frac12\int_{\partial\Omega}\frac{1}{c_\varepsilon}
\frac{\partial}{\partial{\bf n}}|\nabla c_\varepsilon|^2 ds\le \kappa
\int_{\partial\Omega}\frac{1}{c_\varepsilon}|\nabla c_\varepsilon|^2 ds=
\kappa\int_{\partial\Omega}|c_\varepsilon^{\frac12}\nabla\ln c_\varepsilon|^2 ds\nonumber
\\
\le &\eta_1\kappa\int_{\Omega}|D (c_\varepsilon^{\frac12}\nabla\ln c_\varepsilon)|^2 dx+
C_{\eta_1}\int_{\Omega}c_\varepsilon|\nabla\ln c_\varepsilon|^2 dx \nonumber
\\
\le &\eta_1\kappa\int_{\Omega}\left(\frac12 c_\varepsilon^{-\frac12}\nabla c_\varepsilon\nabla\ln c_\varepsilon
+c_\varepsilon^{\frac12} D^2\ln c_\varepsilon\right)^2 dx+
C_{\eta_1}\int_{\Omega}\frac{|\nabla c_\varepsilon|^2}{c_\varepsilon} dx\nonumber
\\
\le & 2\eta_1\kappa\int_{\Omega}
\left(c_\varepsilon | D^2\ln c_\varepsilon|^2+\frac14\frac{|\nabla c_\varepsilon|^4}{c_\varepsilon^3}\right) dx
+C_{\eta_1}\int_{\Omega}\frac{|\nabla c_\varepsilon|^2}{c_\varepsilon} dx\nonumber
\\
\label{3-10}
\le & 10\eta_1\kappa\int_{\Omega} c_\varepsilon | D^2\ln c_\varepsilon|^2 dx
+C_{\eta_1}\int_{\Omega}\frac{|\nabla c_\varepsilon|^2}{c_\varepsilon} dx.
\end{align}
By \eqref{3-3} and Lemma \ref{lem2-5},  for any sufficiently small $\eta_2>0$,  we have that
\begin{align}
\int_{\Omega}\frac{|\nabla c_\varepsilon|^2}{c_\varepsilon} dx &\le
\eta_2\int_{\Omega}\frac{|\nabla c_\varepsilon|^4}{c_\varepsilon^3} dx+C_{\eta_2}\nonumber
\\
\label{3-11}
&\le 16\eta_2\int_{\Omega} c_\varepsilon | D^2\ln c_\varepsilon|^2 dx+C_{\eta_2}.
\end{align}
Substituting \eqref{3-10}, \eqref{3-11}  into \eqref{3-9} gives
\begin{align}
&\frac{d}{dt}\int_\Omega\left(\frac{|\nabla c_\varepsilon|^2}{c_\varepsilon}+
\frac{2}{\chi}n_{\varepsilon}\ln n_{\varepsilon}\right)dx
+\int_\Omega \left(c_{\varepsilon}| D^2\ln c_{\varepsilon}|^2+\frac{n_\varepsilon}{c_\varepsilon}|\nabla c_\varepsilon|^2 +\frac{8}{m\chi}\frac{n_{\varepsilon}+\varepsilon}{n_{\varepsilon}}
|\nabla (n_{\varepsilon}+\varepsilon)^{\frac m2}|^2\right)dx\nonumber
\\
\label{3-12}
& \le \hat C\int_\Omega|\nabla u_\varepsilon|^2 dx+C.
\end{align}
Multiplying the third equation of \eqref{3-1} by $u_\varepsilon$, and integrating it over $\Omega$ yields
\begin{align}
&\frac12\frac{d}{dt}\int_\Omega |u_\varepsilon|^2 dx+\int_\Omega |\nabla u_\varepsilon|^2 dx
=\int_\Omega n_\varepsilon\nabla\varphi_\varepsilon \cdot u_\varepsilon dx
\le \|\nabla\varphi_\varepsilon\|_{L^\infty}\|n_\varepsilon\|_{L^{\frac 65}}\|u_\varepsilon\|_{L^6}\nonumber
\\
&\le C\|n_\varepsilon\|_{L^{\frac 65}}\|\nabla u_\varepsilon\|_{L^2}\nonumber
\\
\label{3-13}
&\le \frac12\int_\Omega |\nabla u_\varepsilon|^2 dx +\tilde C\|n_\varepsilon\|_{L^{\frac 65}}^2.
\end{align}
Using Gagliardo-Nirenberg interpolation inequality, we derive that
\begin{align*}
\tilde C\|n_\varepsilon\|_{L^{\frac 65}}^2& \le\tilde C\|(n_\varepsilon+\varepsilon)^{\frac m2}\|_{L^{\frac{12}{5m}}}^{\frac 4m}
\\
&\le C_1\|\nabla (n_\varepsilon+\varepsilon)^{\frac m2}\|_{L^2}^{\frac{2}{3m-1}}\|(n_\varepsilon+\varepsilon)^{\frac m2}\|_{L^{\frac2m}}^{\frac{10m-4}{m(3m-1)}}+
C_2\|n_\varepsilon+\varepsilon\|_{L^1}^2
\\
&\le C_3\|\nabla (n_\varepsilon+\varepsilon)^{\frac m2}\|_{L^2}^{\frac{2}{3m-1}}+C_4,
\end{align*}
substituting it into \eqref{3-13} yields
\begin{align}\label{3-14}
\frac{d}{dt}\int_\Omega |u_\varepsilon|^2 dx+\int_\Omega |\nabla u_\varepsilon|^2 dx
\le 2C_3\|\nabla (n_\varepsilon+\varepsilon)^{\frac m2}\|_{L^2}^{\frac{2}{3m-1}}+2C_4.
\end{align}
Combining \eqref{3-12} and \eqref{3-14}, and noticing that ${\frac{2}{3m-1}}<2$, then
\begin{align*}
&\frac{d}{dt}\int_\Omega\left(\frac{|\nabla c_\varepsilon|^2}{c_\varepsilon}+
\frac{2}{\chi}n_{\varepsilon}\ln n_{\varepsilon}+2\hat C |u_\varepsilon|^2\right)dx\nonumber
\\
&+\int_\Omega \left(c_{\varepsilon}| D^2\ln c_{\varepsilon}|^2+\frac{n_\varepsilon}{c_\varepsilon}|\nabla c_\varepsilon|^2 +\frac{8}{m\chi}\frac{n_{\varepsilon}+\varepsilon}{n_{\varepsilon}}
|\nabla (n_{\varepsilon}+\varepsilon)^{\frac m2}|^2+\hat C |\nabla u_\varepsilon|^2\right)dx\nonumber
\\
& \le 4\hat C C_3\|\nabla (n_\varepsilon+\varepsilon)^{\frac m2}\|_{L^2}^{\frac{2}{3m-1}}+C_5
\\
&\le \frac{1}{m\chi}\|\nabla (n_\varepsilon+\varepsilon)^{\frac m2}\|_{L^2}^2+C_6,
\end{align*}
that is
\begin{align}
&\frac{d}{dt}\int_\Omega\left(\frac{|\nabla c_\varepsilon|^2}{c_\varepsilon}+
\frac{2}{\chi}n_{\varepsilon}\ln n_{\varepsilon}+2\hat C |u_\varepsilon|^2\right)dx\nonumber
\\
\label{3-15}
&+\int_\Omega \left(c_{\varepsilon}| D^2\ln c_{\varepsilon}|^2+\frac{n_\varepsilon}{c_\varepsilon}|\nabla c_\varepsilon|^2 +\frac{7}{m\chi}\frac{n_{\varepsilon}+\varepsilon}{n_{\varepsilon}}
|\nabla (n_{\varepsilon}+\varepsilon)^{\frac m2}|^2+\hat C |\nabla u_\varepsilon|^2\right)dx\le C.
\end{align}
By Gagliardo-Nirenberg interpolation inequality and \eqref{3-4},  we see that
\begin{align}
\|n_{\varepsilon}\|_{L^{m+\frac23}}^{m+\frac23}\le \|(n_{\varepsilon}+\varepsilon)^{\frac m2}\|_{L^{2+\frac{4}{3m}}}^{2+\frac{4}{3m}}
&\le C_1\|(n_{\varepsilon}+\varepsilon)^{\frac m2}\|_{L^{\frac 2m}}^{\frac{4}{3m}}\|\nabla (n_{\varepsilon}+\varepsilon)^{\frac m2}\|_{L^2}^2
+C_2\|(n_{\varepsilon}+\varepsilon)\|_{L^1}^{m+\frac23}\nonumber
\\
&\le C_1\|(n_{\varepsilon}+\varepsilon)\|_{L^1}^{\frac23}\|\nabla (n_{\varepsilon}+\varepsilon)^{\frac m2}\|_{L^2}^2
+C_2\|(n_{\varepsilon}+\varepsilon)\|_{L^1}^{m+\frac23}\nonumber
\\
\label{3-16}
&=C_3\|\nabla (n_{\varepsilon}+\varepsilon)^{\frac m2}\|_{L^2}^2
+C_4.
\end{align}
Recalling \eqref{3-11}, \eqref{3-16} and Poincar\'e inequality,  there exists $\sigma>0$, such that
\begin{align*}
&\int_\Omega\left(\frac{|\nabla c_\varepsilon|^2}{c_\varepsilon}+
\frac{2}{\chi}n_{\varepsilon}\ln n_{\varepsilon}+2\hat C |u_\varepsilon|^2\right)dx
\\
&\le
\sigma
\int_\Omega \left(c_{\varepsilon}| D^2\ln c_{\varepsilon}|^2+\frac{n_\varepsilon}{c_\varepsilon}|\nabla c_\varepsilon|^2 +\frac{7}{m\chi}
|\nabla (n_{\varepsilon}+\varepsilon)^{\frac m2}|^2+\hat C |\nabla u_\varepsilon|^2\right)dx +C.
\end{align*}
Combining with \eqref{3-15} gives
\begin{align}
&\sup_{t\in(0,T_{\max})}\int_\Omega\left(\frac{|\nabla c_\varepsilon|^2}{c_\varepsilon}+
n_{\varepsilon}\ln n_{\varepsilon}+ |u_\varepsilon|^2\right)dx\nonumber
\nonumber
\\
\label{add3-18}
&+\sup_{t\in(\tau,T_{\max})}\int_{t-\tau}^{t}
\int_\Omega \left(c_{\varepsilon}| D^2\ln c_{\varepsilon}|^2+\frac{n_\varepsilon}{c_\varepsilon}|\nabla c_\varepsilon|^2 +\frac{n_{\varepsilon}+\varepsilon}{n_{\varepsilon}}
|\nabla (n_{\varepsilon}+\varepsilon)^{\frac m2}|^2+ |\nabla u_\varepsilon|^2\right)dx\le C.
\end{align}
Recalling \eqref{add3-7}, we see that
\begin{align*}
&\int_\Omega \frac{| D^2c_{\varepsilon}|^2}{c_{\varepsilon}}dx=\int_\Omega
c_{\varepsilon}| D^2\ln c_{\varepsilon}|^2 dx-\int_\Omega\frac{|\nabla c_\varepsilon|^2}{c_\varepsilon^2}\Delta c_\varepsilon dx
+\int_\Omega\frac{|\nabla c_\varepsilon|^4}{c_\varepsilon^3}dx
\\
&\le \frac14\int_\Omega \frac{|\Delta c_{\varepsilon}|^2}{c_{\varepsilon}}dx+
\int_\Omega
c_{\varepsilon}| D^2\ln c_{\varepsilon}|^2 dx+2\int_\Omega\frac{|\nabla c_\varepsilon|^4}{c_\varepsilon^3}dx
\\
&\le \frac34\int_\Omega \frac{| D^2c_{\varepsilon}|^2}{c_{\varepsilon}}dx+
\int_\Omega
c_{\varepsilon}| D^2\ln c_{\varepsilon}|^2 dx+2\int_\Omega\frac{|\nabla c_\varepsilon|^4}{c_\varepsilon^3}dx.
\end{align*}
Using lemma \ref{lem2-5}, that is
$$
\int_\Omega\frac{|\nabla c_\varepsilon|^4}{c_\varepsilon^3}dx\le C\int_\Omega c_{\varepsilon}| D^2\ln c_{\varepsilon}|^2 dx,
$$
we arrive at
\begin{align*}
\int_\Omega \frac{| D^2c_{\varepsilon}|^2}{c_{\varepsilon}}dx\le C\int_\Omega c_{\varepsilon}| D^2\ln c_{\varepsilon}|^2 dx.
\end{align*}
Combining with \eqref{add3-18}, and this lemma is proved.
\hfill $\Box$

\medskip

Based on this result, we further obtain that

\begin{lemma}
\label{lem3-4}
Assume that $m>1$.
Let $(n_\varepsilon, c_\varepsilon, u_\varepsilon, \pi_\varepsilon)$ be the classical
solution of \eqref{3-1} in $(0, T_{\max})$.  Then we have
\begin{equation}
\label{3-17}
\sup_{t\in(0,+\infty)}\|n_\varepsilon\|_{L^{m}}+\sup_{t\in(0,+\infty)}\int_t^{t+1}
\int_\Omega (n_\varepsilon+\varepsilon)^{2m-3}|\nabla n_\varepsilon|^2 dxds\le C,
\end{equation}
where $C$ is independent of $\varepsilon$.
\end{lemma}

{\it\bf Proof.} Multiplying the first equation of \eqref{3-1}  by $(n_\varepsilon+\varepsilon)^r$ for any $r>0$,
then integrating it over $\Omega$, and using \eqref{3-2},  we obtain
\begin{align}
&\frac1{r+1}\frac{d}{dt}\int_\Omega |n_\varepsilon+\varepsilon|^{r+1} dx
+rm\int_\Omega |n_\varepsilon+\varepsilon|^{m+r-2}|\nabla n_\varepsilon|^2 dx+\int_\Omega |n_\varepsilon+\varepsilon|^{r+1} dx\nonumber
\\
= & r\chi\int_\Omega (n_\varepsilon+\varepsilon)^{r-1}n_\varepsilon\nabla c_\varepsilon
\nabla n_\varepsilon dx+\int_\Omega |n_\varepsilon+\varepsilon|^{r+1} dx\nonumber
\\
\label{3-18}
\le &\int_\Omega |n_\varepsilon+\varepsilon|^{r+1} dx+
\frac{rm}4\int_\Omega |n_\varepsilon+\varepsilon|^{m+r-2}|\nabla n_\varepsilon|^2 dx+
C\int_\Omega |n_\varepsilon+\varepsilon|^{r+2-m}|\nabla c_\varepsilon|^{2}dx.
\end{align}
By Gagliardo-Nirenberg interpolation inequality,  we see that
\begin{align*}
\|n_\varepsilon+\varepsilon\|_{L^{r+1}}^{r+1}&=\|(n_\varepsilon+\varepsilon)^{\frac{m+r}2}\|_{L^{\frac{2(r+1)}{m+r}}}^{\frac{2(r+1)}{m+r}}
\\
&\le C_1\|(n_\varepsilon+\varepsilon)^{\frac{m+r}{2}}\|_{L^{\frac{2}{m+r}}}^{\frac{6m+4r-2}{(m+r)(3m+3r-1)}}
\|\nabla (n_\varepsilon+\varepsilon)^{\frac{m+r}{2}}\|_{L^2}^{\frac{6r}{3m+3r-1}}+C_2\|n_\varepsilon+\varepsilon\|_{L^1}^{r+1}
\\
&\le C_3(1+\|\nabla (n_\varepsilon+\varepsilon)^{\frac{m+r}{2}}\|_{L^2}^{\frac{6r}{3m+3r-1}})
\\
&\le \frac{mr}{(m+r)^2}\|\nabla (n_\varepsilon+\varepsilon)^{\frac{m+r}{2}}\|_{L^2}^2+C_4.
\end{align*}
Substituting the above inequality into \eqref{3-18}, we obtain
\begin{align}
&\frac1{r+1}\frac{d}{dt}\int_\Omega |n_\varepsilon+\varepsilon|^{r+1} dx
+\frac{rm}2\int_\Omega |n_\varepsilon+\varepsilon|^{m+r-2}|\nabla n_\varepsilon|^2 dx+\int_\Omega |n_\varepsilon+\varepsilon|^{r+1} dx
\nonumber
\\
\label{3-19}
\le & C\int_\Omega |n_\varepsilon+\varepsilon|^{r+2-m}|\nabla c_\varepsilon|^{2}dx+C_4.
\end{align}
Taking $r=m-1$ in \eqref{3-19}, and noticing
that $$
(n_\varepsilon+\varepsilon)|\nabla c_\varepsilon|^{2}\le \|c_\varepsilon\|_{L^\infty}
\frac{n_\varepsilon}{c_\varepsilon}|\nabla c_\varepsilon|^{2}+\varepsilon\|c_\varepsilon\|_{L^\infty}\frac{|\nabla c_\varepsilon|^{2}}{c_\varepsilon},
$$
then combining with  \eqref{3-5}, it gives \eqref{3-17}.  \hfill $\Box$

\medskip

In order to avoid the tedious discussion of different situations, in what follows, we may assume that $1<m\le\frac 54$ since  that the global  existence had been proved for the case $m>\frac54$ in some previous paper \cite{W2, W3}.

\begin{lemma}
\label{lem3-5}
Assume  $1<m\le \frac 54$.
Let $(n_\varepsilon, c_\varepsilon, u_\varepsilon, \pi_\varepsilon)$ be the classical
solution of \eqref{3-1} in $(0, T_{\max})$. Then
\begin{align}
&\frac{d}{dt}\int_\Omega\frac{|\nabla n_\varepsilon|^2}{n_\varepsilon+\varepsilon}dx
+\frac{2m}{(m-1)^2}\int_\Omega (n_\varepsilon+\varepsilon)| D^2 (n_\varepsilon+\varepsilon)^{\frac{m-1}2}|^2 dx\nonumber
\\
&\le C \int_{\Omega}(n_\varepsilon+\varepsilon)^{m-2}|\nabla n_\varepsilon |^2 dx+C\int_\Omega(n_\varepsilon+\varepsilon)^{4-3m} |u_{\varepsilon}|^4dx\nonumber
\\
\label{3-20}
&+C\int_\Omega\left((n_\varepsilon+\varepsilon)^{4-3m} |\nabla c_{\varepsilon}|^4+(n_\varepsilon+\varepsilon)^{2-m}|\Delta c_\varepsilon|^2\right) dx
\end{align}
for  $t\in(0, T_{\max})$,
where $C$ depends on $r$, and  which is independent of $\varepsilon$.
\end{lemma}

{\it\bfseries Proof. } By a direct calculation, we obtain that
\begin{align}
&\frac{d}{dt}\int_\Omega\frac{|\nabla n_\varepsilon|^2}{n_\varepsilon+\varepsilon}dx
=2\int_\Omega\frac{\nabla n_\varepsilon}{n_\varepsilon+\varepsilon}\nabla n_{\varepsilon t}dx
-\int_\Omega\frac{|\nabla n_\varepsilon|^2}{(n_\varepsilon+\varepsilon)^2}n_{\varepsilon t}dx
\nonumber
\\
&=-2\int_\Omega n_{\varepsilon t}\left(\frac{\Delta n_{\varepsilon}}{n_\varepsilon+\varepsilon}
-\frac{|\nabla n_\varepsilon|^2}{(n_\varepsilon+\varepsilon)^2}\right)dx
-\int_\Omega\frac{|\nabla n_\varepsilon|^2}{(n_\varepsilon+\varepsilon)^2}n_{\varepsilon t}dx
\nonumber
\\
&=\int_\Omega \left(-2\frac{\Delta n_{\varepsilon}}{n_\varepsilon+\varepsilon}+\frac{|\nabla n_\varepsilon|^2}{(n_\varepsilon+\varepsilon)^2}\right) n_{\varepsilon t}dx\nonumber
\\
&=\int_\Omega \left(-2\frac{\Delta n_{\varepsilon}}{n_\varepsilon+\varepsilon}+\frac{|\nabla n_\varepsilon|^2}{(n_\varepsilon+\varepsilon)^2}\right) (m\nabla\cdot((\varepsilon +n_{\varepsilon})^{m-1}\nabla n_{\varepsilon} )-u_{\varepsilon}\cdot\nabla n_{\varepsilon}-\chi
\nabla\cdot\left(n_{\varepsilon}\nabla c_{\varepsilon}\right)) dx\nonumber
\\
\label{3-21}
&=I+II+III,
\end{align}
where
\begin{align*}
&I=m\int_\Omega \left(-2\frac{\Delta n_{\varepsilon}}{n_\varepsilon+\varepsilon}+\frac{|\nabla n_\varepsilon|^2}{(n_\varepsilon+\varepsilon)^2}\right) \nabla\cdot((\varepsilon +n_{\varepsilon})^{m-1}\nabla n_{\varepsilon} ) dx,
\\
&II=\int_\Omega \left(2\frac{\Delta n_{\varepsilon}}{n_\varepsilon+\varepsilon}-\frac{|\nabla n_\varepsilon|^2}{(n_\varepsilon+\varepsilon)^2}\right)
u_{\varepsilon}\cdot\nabla n_{\varepsilon} dx,
\\
&III=\chi\int_\Omega \left(2\frac{\Delta n_{\varepsilon}}{n_\varepsilon+\varepsilon}-\frac{|\nabla n_\varepsilon|^2}{(n_\varepsilon+\varepsilon)^2}\right)
\nabla\cdot(n_{\varepsilon}\nabla c_{\varepsilon})  dx.
\end{align*}
Next, we estimate the three terms respectively.
Firstly, for $I$, it is easy to see that
\begin{align}
I=&-2m\int_\Omega (n_\varepsilon+\varepsilon)^{m-2}|\Delta n_{\varepsilon}|^2 dx+(3m-2m^2)\int_\Omega (n_\varepsilon+\varepsilon)^{m-3}|\nabla n_{\varepsilon} |^2\Delta n_{\varepsilon}dx\nonumber
\\
&+m(m-1)\int_\Omega (n_\varepsilon+\varepsilon)^{m-4}|\nabla n_{\varepsilon} |^4dx. \label{3-22}
\end{align}
By Lemma \ref{lem2-4},
\begin{align}
-2m\int_\Omega (n_\varepsilon+\varepsilon)^{m-2}|\Delta n_{\varepsilon}|^2 dx=3m(m-2)\int_\Omega (n_\varepsilon+\varepsilon)^{m-3}|\nabla n_{\varepsilon}|^2\Delta n_{\varepsilon} dx-2m\int_\Omega(n_\varepsilon+\varepsilon)^{m-2}| D^2 n_{\varepsilon}|^2dx\nonumber
\\
\label{3-23}
+m(m-2)(m-3)
\int_\Omega(n_\varepsilon+\varepsilon)^{m-4}|\nabla n_{\varepsilon}|^4dx+m\int_{\partial\Omega}
(n_\varepsilon+\varepsilon)^{m-2}\frac{\partial|\nabla n_{\varepsilon}|^2}{\partial\nu}dS.
\end{align}
From a direct calculation, we infer that
\allowdisplaybreaks
\begin{align*}
&\int_\Omega (n_\varepsilon+\varepsilon)| D^2 (n_\varepsilon+\varepsilon)^{\frac{m-1}2}|^2 dx
\\
&=\int_\Omega (n_\varepsilon+\varepsilon)\left(\frac{m-1}2(n_\varepsilon+\varepsilon)^{\frac{m-3}2} D^2 n_\varepsilon+\frac{(m-1)(m-3)}4(n_\varepsilon+\varepsilon)^{\frac{m-5}2}\nabla n_\varepsilon\otimes\nabla n_\varepsilon\right)^2dx
\\
&=\frac{(m-1)^2}4\int_\Omega \left((n_\varepsilon+\varepsilon)^{m-2}| D^2 n_\varepsilon|^2
+\frac{(m-3)^2}{4}(n_\varepsilon+\varepsilon)^{m-4}|\nabla n_\varepsilon|^4+(m-3)(n_\varepsilon+\varepsilon)^{m-3}\nabla n_\varepsilon\cdot D^2 n_\varepsilon\cdot\nabla n_\varepsilon\right)dx
\\
&=\frac{(m-1)^2}4\int_\Omega \left((n_\varepsilon+\varepsilon)^{m-2}| D^2 n_\varepsilon|^2
+\frac{(m-3)^2}{4}(n_\varepsilon+\varepsilon)^{m-4}|\nabla n_\varepsilon|^4-\frac{m-3}2|\nabla n_\varepsilon|^2 \nabla\cdot((n_\varepsilon+\varepsilon)^{m-3}\nabla n_\varepsilon)\right)dx
\\
&=\frac{(m-1)^2}4\int_\Omega \left((n_\varepsilon+\varepsilon)^{m-2}| D^2 n_\varepsilon|^2
-\frac{(m-3)^2}{4}(n_\varepsilon+\varepsilon)^{m-4}|\nabla n_\varepsilon|^4-\frac{m-3}2|\nabla n_\varepsilon|^2 (n_\varepsilon+\varepsilon)^{m-3}\Delta n_\varepsilon\right)dx,
\end{align*}
where $|D^2 n_\varepsilon|^2=\sum_{i,j=1}^n |D_{ij}n_\varepsilon|^2$, $\nabla n_\varepsilon\cdot D^2 n_\varepsilon\cdot\nabla n_\varepsilon=\sum_{i,j=1}^n D_i n_\varepsilon D_{ij}n_\varepsilon D_jn_\varepsilon$,
that is
\begin{align}
&\int_\Omega (n_\varepsilon+\varepsilon)^{m-2}| D^2 n_\varepsilon|^2 dx=\frac 4{(m-1)^2}\int_\Omega (n_\varepsilon+\varepsilon)| D^2 (n_\varepsilon+\varepsilon)^{\frac{m-1}2}|^2 dx\nonumber
\\
\label{3-24}
&+\frac{(m-3)^2}{4}\int_\Omega(n_\varepsilon+\varepsilon)^{m-4}|\nabla n_\varepsilon|^4dx+
\frac{m-3}2\int_\Omega|\nabla n_\varepsilon|^2 (n_\varepsilon+\varepsilon)^{m-3}\Delta n_\varepsilon dx.
\end{align}
Substituting \eqref{3-24} into \eqref{3-23}, then substituting the result into \eqref{3-22} again, and using Lemma \ref{lem2-6}, we get that
\begin{align}
&I=-\frac{8m}{(m-1)^2}\int_\Omega (n_\varepsilon+\varepsilon)| D^2 (n_\varepsilon+\varepsilon)^{\frac{m-1}2}|^2 dx
+\frac{m(m-1)^2}2\int_\Omega(n_\varepsilon+\varepsilon)^{m-4}|\nabla n_\varepsilon|^4dx\nonumber
\\
&+m\int_{\partial\Omega}
(n_\varepsilon+\varepsilon)^{m-2}\frac{\partial|\nabla n_{\varepsilon}|^2}{\partial\nu}dS\nonumber
\\
&\le -\frac{8m}{(m-1)^2}\int_\Omega (n_\varepsilon+\varepsilon)| D^2 (n_\varepsilon+\varepsilon)^{\frac{m-1}2}|^2 dx
+\frac{m(m-1)^2}2\int_\Omega(n_\varepsilon+\varepsilon)^{m-4}|\nabla n_\varepsilon|^4dx \nonumber
\\
&+2\kappa m\int_{\partial\Omega}
(n_\varepsilon+\varepsilon)^{m-2}|\nabla n_{\varepsilon}|^2dS\nonumber
\\
&\le -\frac{8m}{(m-1)^2}\int_\Omega (n_\varepsilon+\varepsilon)| D^2 (n_\varepsilon+\varepsilon)^{\frac{m-1}2}|^2 dx
+\frac{m(m-1)^2}2\int_\Omega(n_\varepsilon+\varepsilon)^{m-4}|\nabla n_\varepsilon|^4dx \nonumber
\\
\label{3-25}
&+\frac{8\kappa m}{(m-1)^2}\int_{\partial\Omega}
(n_\varepsilon+\varepsilon)|\nabla (n_\varepsilon+\varepsilon)^{\frac{m-1}{2}}|^2dS.
\end{align}
Recalling Lemma \ref{lem2-5} with $h(\varphi)=\varphi^{\frac{3-m}2}$, then
we have
\begin{equation}\label{3-26}
\int_\Omega(n_\varepsilon+\varepsilon)^{m-4}|\nabla n_\varepsilon|^4dx\le \frac{16(2+\sqrt 3)^2}{(3-m)^2(m-1)^2}
\int_\Omega (n_\varepsilon+\varepsilon)| D^2 (n_\varepsilon+\varepsilon)^{\frac{m-1}2}|^2 dx.
\end{equation}
Noticing that when $m\le \frac 54$,
$$
\frac{8m(2+\sqrt 3)^2}{(3-m)^2}<\frac{4m}{(m-1)^2},
$$
then
\begin{align}
&\frac{m(m-1)^2}2\int_\Omega(n_\varepsilon+\varepsilon)^{m-4}|\nabla n_\varepsilon|^4dx\nonumber
\\
&\le \frac{8m(2+\sqrt 3)^2}{(3-m)^2}
\int_\Omega (n_\varepsilon+\varepsilon)| D^2 (n_\varepsilon+\varepsilon)^{\frac{m-1}2}|^2 dx\nonumber
\\
\label{3-27}
&\le \frac{4m}{(m-1)^2}\int_\Omega (n_\varepsilon+\varepsilon)| D^2 (n_\varepsilon+\varepsilon)^{\frac{m-1}2}|^2 dx.
\end{align}
By boundary trace imbedding theorem, for any small $\eta>0$, we obtain that
\begin{align}
&\int_{\partial\Omega}
(n_\varepsilon+\varepsilon)|\nabla (n_\varepsilon+\varepsilon)^{\frac{m-1}{2}}|^2dS\nonumber
\\
&\le \eta\int_{\Omega}\left|\nabla\left((n_\varepsilon+\varepsilon)^{\frac 12}\nabla (n_\varepsilon+\varepsilon)^{\frac{m-1}{2}}\right) \right|^2 dx
+C_\eta\int_{\Omega}(n_\varepsilon+\varepsilon)|\nabla (n_\varepsilon+\varepsilon)^{\frac{m-1}{2}}|^2 dx \nonumber
\\
&\le \eta\int_{\Omega}\left|(n_\varepsilon+\varepsilon)^{\frac 12} D^2 (n_\varepsilon+\varepsilon)^{\frac{m-1}{2}}+\frac{m-1}4
(n_\varepsilon+\varepsilon)^{\frac m2-2}|\nabla n_\varepsilon|^2\right|^2 dx
+\frac{C_\eta(m-1)^2}4 \int_{\Omega}(n_\varepsilon+\varepsilon)^{m-2}|\nabla n_\varepsilon |^2 dx \nonumber
\\
&\le 2\eta\int_{\Omega}\left((n_\varepsilon+\varepsilon)\left| D^2 (n_\varepsilon+\varepsilon)^{\frac{m-1}{2}}\right|^2+\frac{(m-1)^2}{16}
(n_\varepsilon+\varepsilon)^{m-4}|\nabla n_\varepsilon|^4 \right)dx \nonumber
\\
&\quad +\frac{C_\eta(m-1)^2}4 \int_{\Omega}(n_\varepsilon+\varepsilon)^{m-2}|\nabla n_\varepsilon |^2 dx\nonumber
\\
\label{3-28}
&\le  2\eta\left(1+\frac{(2+\sqrt 3)^2}{(3-m)^2}\right)\int_{\Omega}(n_\varepsilon+\varepsilon)\left| D^2 (n_\varepsilon+\varepsilon)^{\frac{m-1}{2}}\right|^2dx
+\frac{C_\eta(m-1)^2}4 \int_{\Omega}(n_\varepsilon+\varepsilon)^{m-2}|\nabla n_\varepsilon |^2 dx.
\end{align}
Combining \eqref{3-25}, \eqref{3-27} and \eqref{3-28}, and taking $\eta$ appropriately small, we finally obtain that
\begin{align}\label{3-29}
I\le -\frac{3m}{(m-1)^2}\int_\Omega (n_\varepsilon+\varepsilon)| D^2 (n_\varepsilon+\varepsilon)^{\frac{m-1}2}|^2 dx
+C \int_{\Omega}(n_\varepsilon+\varepsilon)^{m-2}|\nabla n_\varepsilon |^2 dx.
\end{align}
Next, we turn our attention to $II$. For this purpose, recalling \eqref{3-24}, \eqref{3-26}, and noticing that $|\Delta n_\varepsilon|^2\le 3| D^2 n_\varepsilon|^2$, then
\begin{align*}
&\int_\Omega (n_\varepsilon+\varepsilon)^{m-2}|\Delta n_\varepsilon|^2 dx\le 3\int_\Omega (n_\varepsilon+\varepsilon)^{m-2}| D^2 n_\varepsilon|^2 dx
\\
&\le \frac{12}{(m-1)^2}\int_\Omega (n_\varepsilon+\varepsilon)| D^2 (n_\varepsilon+\varepsilon)^{\frac{m-1}2}|^2 dx+\frac{3(m-3)^2}{4}\int_\Omega(n_\varepsilon+\varepsilon)^{m-4}|\nabla n_\varepsilon|^4dx
\\
&+\frac12\int_\Omega (n_\varepsilon+\varepsilon)^{m-2}|\Delta n_\varepsilon|^2 dx+
\frac{9(m-3)^2}8\int_\Omega|\nabla n_\varepsilon|^4 (n_\varepsilon+\varepsilon)^{m-4}dx
\\
&\le \left(\frac{12}{(m-1)^2}+\frac{30(2+\sqrt 3)^2}{(m-1)^2}\right)\int_\Omega (n_\varepsilon+\varepsilon)| D^2 (n_\varepsilon+\varepsilon)^{\frac{m-1}2}|^2 dx+\frac12\int_\Omega (n_\varepsilon+\varepsilon)^{m-2}|\Delta n_\varepsilon|^2 dx,
\end{align*}
which implies that
\begin{align}\label{3-30}
&\int_\Omega (n_\varepsilon+\varepsilon)^{m-2}|\Delta n_\varepsilon|^2 dx
\le \left(\frac{24}{(m-1)^2}+\frac{60(2+\sqrt 3)^2}{(m-1)^2}\right)\int_\Omega (n_\varepsilon+\varepsilon)| D^2 (n_\varepsilon+\varepsilon)^{\frac{m-1}2}|^2 dx.
\end{align}
For $II$, for any small $\eta>0$, we have
\begin{align*}
&II=\int_\Omega \left(2\frac{\Delta n_{\varepsilon}}{n_\varepsilon+\varepsilon}-\frac{|\nabla n_\varepsilon|^2}{(n_\varepsilon+\varepsilon)^2}\right)
u_{\varepsilon}\cdot\nabla n_{\varepsilon} dx
\\
&\le \eta\int_\Omega (n_\varepsilon+\varepsilon)^{m-2}|\Delta n_\varepsilon|^2 dx+C_{\eta}\int_\Omega(n_\varepsilon+\varepsilon)^{-m} |u_{\varepsilon}|^2|\nabla n_{\varepsilon}|^2dx
+\frac{\eta}2\int_\Omega (n_\varepsilon+\varepsilon)^{m-4}|\nabla n_\varepsilon|^4 dx
\\
&\le  \eta\int_\Omega (n_\varepsilon+\varepsilon)^{m-2}|\Delta n_\varepsilon|^2 dx
+\eta\int_\Omega (n_\varepsilon+\varepsilon)^{m-4}|\nabla n_\varepsilon|^4 dx+\tilde C_\eta
\int_\Omega(n_\varepsilon+\varepsilon)^{4-3m} |u_{\varepsilon}|^4dx.
\end{align*}
Combining \eqref{3-26}, \eqref{3-30}, and taking $\eta$ appropriately small in the above inequality, we infer that
\begin{align}\label{3-31}
&II\le \frac{m}{2(m-1)^2}\int_\Omega (n_\varepsilon+\varepsilon)| D^2 (n_\varepsilon+\varepsilon)^{\frac{m-1}2}|^2 dx
+C\int_\Omega(n_\varepsilon+\varepsilon)^{4-3m} |u_{\varepsilon}|^4dx.
\end{align}
For $III$, we see that
$$
III=III_1+III_2,
$$
where
$$
III_1=\chi\int_\Omega \left(2\frac{\Delta n_{\varepsilon}}{n_\varepsilon+\varepsilon}-\frac{|\nabla n_\varepsilon|^2}{(n_\varepsilon+\varepsilon)^2}\right)
\nabla n_{\varepsilon}\nabla c_{\varepsilon}dx,
$$
$$
III_2=\chi\int_\Omega \left(2\frac{\Delta n_{\varepsilon}}{n_\varepsilon+\varepsilon}-\frac{|\nabla n_\varepsilon|^2}{(n_\varepsilon+\varepsilon)^2}\right)
n_{\varepsilon}\Delta c_{\varepsilon}dx.
$$
For $III_1$, it completely similar to the estimation for $II$, then for any small constant $\eta>0$,
\begin{align*}
&III \le  \eta\int_\Omega (n_\varepsilon+\varepsilon)^{m-2}|\Delta n_\varepsilon|^2 dx
+\eta\int_\Omega (n_\varepsilon+\varepsilon)^{m-4}|\nabla n_\varepsilon|^4 dx+C_\eta
\int_\Omega(n_\varepsilon+\varepsilon)^{4-3m} |\nabla c_{\varepsilon}|^4dx+III_2
\\
&\le 2\eta\int_\Omega (n_\varepsilon+\varepsilon)^{m-2}|\Delta n_\varepsilon|^2 dx
+2\eta\int_\Omega (n_\varepsilon+\varepsilon)^{m-4}|\nabla n_\varepsilon|^4 dx+C_\eta
\int_\Omega(n_\varepsilon+\varepsilon)^{4-3m} |\nabla c_{\varepsilon}|^4dx
\\
&+C_\eta
\int_\Omega (n_\varepsilon+\varepsilon)^{2-m}|\Delta c_\varepsilon|^2 dx
\end{align*}
Using \eqref{3-26}, \eqref{3-30}, and taking $\eta$ appropriately small in the above inequality, it yields
\begin{align}
III\le & \frac{m}{2(m-1)^2}\int_\Omega (n_\varepsilon+\varepsilon)| D^2 (n_\varepsilon+\varepsilon)^{\frac{m-1}2}|^2 dx\nonumber
\\
\label{3-32}
&+C\int_\Omega\left((n_\varepsilon+\varepsilon)^{4-3m} |\nabla c_{\varepsilon}|^4+(n_\varepsilon+\varepsilon)^{2-m}|\Delta c_\varepsilon|^2\right) dx.
\end{align}
Combining \eqref{3-21}, \eqref{3-29}, \eqref{3-31} and \eqref{3-32}, and this lemma is proved. \hfill $\Box$

\begin{lemma}
\label{lem3-6}
Assume  $1<m\le \frac54$.
Let $(n_\varepsilon, c_\varepsilon, u_\varepsilon, \pi_\varepsilon)$ be the classical
solution of \eqref{3-1} in $(0, T_{\max})$. Then
\begin{align}
\label{3-33}
\sup_{t\in(\tau, T_{\max})}\int_{t-\tau}^{t}\int_\Omega (n_{\varepsilon}+\varepsilon)^{2m-4}|\nabla n_{\varepsilon}|^2dxds
\le C,
\end{align}
where $C$  is independent of $\varepsilon$ and $T_{\max}$.
\end{lemma}

 {\bfseries Proof.} Multiplying the first equation of \eqref{3-1} by $-(n_{\varepsilon}+\varepsilon)^{m-2}$,  integrating it over $\Omega$,
 and using \eqref{3-5} yields
\begin{align*}
&-\frac1{m-1}\frac{d}{dt}\int_\Omega (n_{\varepsilon}+\varepsilon)^{m-1}dx+m(2-m)\int_\Omega (n_{\varepsilon}+\varepsilon)^{2m-4}|\nabla n_{\varepsilon}|^2dx
\\
&=\chi(2-m)\int_\Omega n_{\varepsilon}(n_{\varepsilon}+\varepsilon)^{m-3}
\nabla c_{\varepsilon}\nabla n_{\varepsilon} dx
\\
&\le \frac{m(2-m)}2\int_\Omega (n_{\varepsilon}+\varepsilon)^{2m-4}|\nabla n_{\varepsilon}|^2dx+\frac{\chi^2}{2m}\|c_{\varepsilon}\|_{L^\infty}
\int_\Omega \frac{|\nabla c_{\varepsilon}|^2}{c_{\varepsilon}}dx
\\
&\le  \frac{m(2-m)}2\int_\Omega (n_{\varepsilon}+\varepsilon)^{2m-4}|\nabla n_{\varepsilon}|^2dx+C.
\end{align*}
Noticing $0<m-1<1$, using \eqref{3-4}, we complete the proof by a direct integration.  \hfill $\Box$

\begin{lemma}
\label{lem3-7}
Assume that $1<m\le \frac54$.
Let $(n_\varepsilon, c_\varepsilon, u_\varepsilon, \pi_\varepsilon)$ be the classical
solution of \eqref{3-1} in $(0, T_{\max})$. Then for any small constant $\tilde\eta>0$,
\begin{align}
&\frac{d}{dt}\int_\Omega |\nabla c_\varepsilon|^{4} dx
+\int_{\Omega}|\nabla c_\varepsilon|^{2}| D^2 c_\varepsilon|^2 dx+2
\int_{\Omega}|\nabla c_\varepsilon|^{2}(\nabla|\nabla c_\varepsilon|)^2 dx
+2\int_{\Omega} n_\varepsilon|\nabla c_\varepsilon|^{4}dx\nonumber
\\
\label{3-34}
&\le \tilde\eta\int_{\Omega}(n_\varepsilon+\varepsilon)^{m-4}|\nabla n_\varepsilon|^{4}dx+
C_{\tilde\eta}\int_{\Omega}|\nabla c_\varepsilon|^{4}dx
+C\int_{\Omega}|u_\varepsilon|^6 dx
+C,
\end{align}
where $C$, $C_{\tilde\eta}$  are independent of $\varepsilon$ and $T_{\max}$, and $C_{\tilde\eta}$ depends on $\tilde\eta$.
\end{lemma}

{\bfseries Proof.}
Applying $\nabla$ to the second equation of \eqref{3-1}, and
multiplying the resulting equation by $|\nabla c_\varepsilon|^{2}\nabla c_\varepsilon$, and using Lemma \ref{lem2-6}, we obtain
\begin{align}
&\frac14\frac{d}{dt}\int_\Omega |\nabla c_\varepsilon|^{4} dx
+\int_{\Omega}|\nabla c_\varepsilon|^{2}| D^2 c_\varepsilon|^2 dx+2
\int_{\Omega}|\nabla c_\varepsilon|^{2}(\nabla|\nabla c_\varepsilon|)^2 dx
+\int_{\Omega} n_\varepsilon|\nabla c_\varepsilon|^{4}dx\nonumber
\\
&=-\int_{\Omega}c_\varepsilon |\nabla c_\varepsilon|^{2}\nabla c_\varepsilon\nabla n_\varepsilon dx
+\int_{\Omega}u_\varepsilon\nabla c_\varepsilon {\rm div}(|\nabla c_\varepsilon|^{2}\nabla c_\varepsilon)dx+
\frac12\int_{\partial\Omega}\frac{\partial(|\nabla c_\varepsilon|^2)}{\partial{\bf n}}|\nabla c_\varepsilon|^{2} dS\nonumber
\\
&\le \tilde\eta\int_{\Omega}(n_\varepsilon+\varepsilon)^{m-4}|\nabla n_\varepsilon|^{4}dx+
C_{\tilde\eta}\int_{\Omega}(n_\varepsilon+\varepsilon)^{\frac{4-m}3}|\nabla c_\varepsilon|^{4}dx
+\eta\int_{\Omega}|\nabla c_\varepsilon|^6 dx+C_\eta\int_{\Omega}|u_\varepsilon|^6 dx\nonumber
\\
\label{3-35}
&+\frac12\int_{\Omega}|\nabla c_\varepsilon|^{2}| D^2 c_\varepsilon|^2 dx
+\frac 12\int_{\Omega}|\nabla c_\varepsilon|^{2}(\nabla|\nabla c_\varepsilon|)^2 dx+
\kappa\int_{\partial\Omega}|\nabla c_\varepsilon|^4 dS
\end{align}
for any small $\eta, \tilde\eta>0$, and $C_\eta$ depends on $\eta$, $C_{\tilde\eta}$ depends on $\tilde\eta$.
By the boundary trace embedding inequalities, and combining with \eqref{3-5},
we conclude that for any small $\eta>0$,
\begin{align*}
\kappa\int_{\partial\Omega}|\nabla c_\varepsilon|^4 dS
\le \eta\|\nabla(|\nabla c_\varepsilon|^{2})\|_{L^2}^2
+C_\eta\||\nabla c_\varepsilon|^{2}\|_{L^{1}}^2
\\
\le \eta\|\nabla(|\nabla c_\varepsilon|^{2})\|_{L^2}^2
+\hat C_\eta.
\end{align*}
Substituting this inequality into \eqref{3-22} with $\eta$ appropriately small yields
\begin{align}
&\frac14\frac{d}{dt}\int_\Omega |\nabla c_\varepsilon|^{4} dx
+\frac12\int_{\Omega}|\nabla c_\varepsilon|^{2}| D^2 c_\varepsilon|^2 dx+
\int_{\Omega}|\nabla c_\varepsilon|^{2}(\nabla|\nabla c_\varepsilon|)^2 dx
+\int_{\Omega} n_\varepsilon|\nabla c_\varepsilon|^{4}dx\nonumber
\\
\label{3-36}
&\le \tilde\eta\int_{\Omega}(n_\varepsilon+\varepsilon)^{m-4}|\nabla n_\varepsilon|^{4}dx+
C_{\tilde\eta}\int_{\Omega}(n_\varepsilon+\varepsilon)^{\frac{4-m}3}|\nabla c_\varepsilon|^{4}dx
+\eta\int_{\Omega}|\nabla c_\varepsilon|^6 dx+C_\eta\int_{\Omega}|u_\varepsilon|^6 dx
+C.
\end{align}
Noticing that
\begin{align*}
\|\nabla c_\varepsilon\|_{L^6}^{6}
=\int_\Omega |\nabla c_\varepsilon|^{4}\nabla c_\varepsilon\nabla c_\varepsilon dx
=-\int_\Omega c_\varepsilon\left(|\nabla c_\varepsilon|^{4}\Delta c_\varepsilon+
4|\nabla c_\varepsilon|^{3}\nabla c_\varepsilon\nabla|\nabla c_\varepsilon| \right) dx
\\
\le C\left(\int_\Omega |\nabla c_\varepsilon|^{2}| D^2 c_\varepsilon|^2 dx\right)^{\frac12}
\left(\int_\Omega |\nabla c_\varepsilon|^{6}dx\right)^{\frac12}+
\left(\int_{\Omega}|\nabla c_\varepsilon|^{2}(\nabla|\nabla c_\varepsilon|)^2 dx\right)^{\frac12}
\left(\int_\Omega |\nabla c_\varepsilon|^{6}dx\right)^{\frac12},
\end{align*}
then we have
\begin{equation}
\label{3-37}
\|\nabla c_\varepsilon\|_{L^6}^{6}\le
C\int_\Omega |\nabla c_\varepsilon|^{2}| D^2 c_\varepsilon|^2 dx
+C\int_{\Omega}|\nabla c_\varepsilon|^{2}(\nabla|\nabla c_\varepsilon|)^2 dx.
\end{equation}
Combining with \eqref{3-36}, and taking $\eta$ appropriately small, we arrive at
\begin{align*}
&\frac14\frac{d}{dt}\int_\Omega |\nabla c_\varepsilon|^{4} dx
+\frac14\int_{\Omega}|\nabla c_\varepsilon|^{2}| D^2 c_\varepsilon|^2 dx+\frac12
\int_{\Omega}|\nabla c_\varepsilon|^{2}(\nabla|\nabla c_\varepsilon|)^2 dx
+\int_{\Omega} n_\varepsilon|\nabla c_\varepsilon|^{4}dx
\\
&\le \tilde\eta\int_{\Omega}(n_\varepsilon+\varepsilon)^{m-4}|\nabla n_\varepsilon|^{4}dx+
C_{\tilde\eta}\int_{\Omega}(n_\varepsilon+\varepsilon)^{\frac{4-m}3}|\nabla c_\varepsilon|^{4}dx
+C\int_{\Omega}|u_\varepsilon|^6 dx
+C
\\
&\le \tilde\eta\int_{\Omega}(n_\varepsilon+\varepsilon)^{m-4}|\nabla n_\varepsilon|^{4}dx+\frac12
\int_{\Omega}(n_\varepsilon+\varepsilon)|\nabla c_\varepsilon|^{4}dx+
\tilde C_{\tilde\eta}\int_{\Omega}|\nabla c_\varepsilon|^{4}dx
+C\int_{\Omega}|u_\varepsilon|^6 dx
+C,
\end{align*}
since $0<\frac{4-m}3<1$,
and this lemma is proved. \hfill $\Box$

\begin{lemma}
\label{lem3-8}
Assume that $1<m\le \frac54$.
Let $(n_\varepsilon, c_\varepsilon, u_\varepsilon, \pi_\varepsilon)$ be the classical
solution of \eqref{3-1} in $(0, T_{\max})$. Then for any small constant $\eta>0$,
\begin{align}
&\frac{d}{dt}\int_\Omega \left(|\Delta c_\varepsilon|^{2}+|\nabla c_\varepsilon|^{4}\right) dx
+\frac12\int_{\Omega}|\nabla\Delta c_\varepsilon|^{2} dx+\int_{\Omega}n_\varepsilon |\Delta c_\varepsilon|^{2}dx
+\int_{\Omega} n_\varepsilon|\nabla c_\varepsilon|^{4}dx\nonumber
\\
&+\frac12\int_{\Omega}|\nabla c_\varepsilon|^{2}| D^2 c_\varepsilon|^2 dx+
\int_{\Omega}|\nabla c_\varepsilon|^{2}(\nabla|\nabla c_\varepsilon|)^2 dx+\int_\Omega \left(|\Delta c_\varepsilon|^{2}+|\nabla c_\varepsilon|^{4}\right) dx \nonumber
\\
&\le\eta\int_{\Omega}\left((n_\varepsilon+\varepsilon)^{m-2} |\Delta n_\varepsilon|^{2}+(n_\varepsilon+\varepsilon)^{m-4}|\nabla n_\varepsilon|^{4}\right)dx +C_\eta \int_{\Omega}\left(|\Delta c_\varepsilon|^{2}+ |\nabla c_\varepsilon|^{4} \right)dx\nonumber
\\
\label{3-38}
&+C\int_{\Omega}\left(|\nabla u_\varepsilon|^3+|u_\varepsilon|^6\right)dx+C,
\end{align}
where  $C$, $C_\eta$ are independent of $\varepsilon$ and $T_{\max}$, and $C_\eta$ depends on $\eta$.
\end{lemma}

{\bfseries Proof.}
Applying $\nabla$ to the second equation of \eqref{3-1},
multiplying the resulting equation by $-\nabla\Delta c_\varepsilon$ on both sides, and integrating it over $\Omega$ yields
\allowdisplaybreaks
\begin{align*}
&\frac12\frac{d}{dt}\int_\Omega |\Delta c_\varepsilon|^{2} dx
+\int_{\Omega}|\nabla\Delta c_\varepsilon|^{2} dx
\\
&=\int_{\Omega}\nabla(u_\varepsilon\cdot\nabla c_\varepsilon)\nabla\Delta c_\varepsilon dx
-\int_{\Omega}\Delta(n_\varepsilon c_\varepsilon)\Delta c_\varepsilon dx
\\
&= \int_{\Omega}(\nabla u_\varepsilon \nabla c_\varepsilon+u_\varepsilon  D^2c_\varepsilon )\nabla\Delta c_\varepsilon dx
-\int_{\Omega}(c_\varepsilon\Delta n_\varepsilon\Delta c_\varepsilon+n_\varepsilon|\Delta c_\varepsilon|^2
+2\nabla n_\varepsilon\nabla c_\varepsilon\Delta c_\varepsilon) dx
\\
&\le -\int_{\Omega}n_\varepsilon |\Delta c_\varepsilon|^{2}dx
+\frac12\int_{\Omega}|\nabla\Delta c_\varepsilon|^{2} dx+ \int_{\Omega}(|\nabla u_\varepsilon \nabla c_\varepsilon|^2
+|u_\varepsilon  D^2c_\varepsilon |^2) dx+\eta\int_{\Omega}(n_\varepsilon+\varepsilon)^{m-2} |\Delta n_\varepsilon|^{2}dx
\\
&+C_\eta \int_{\Omega}(n_\varepsilon+\varepsilon)^{2-m} |\Delta c_\varepsilon|^{2} dx+\frac14\int_{\Omega}(n_\varepsilon+\varepsilon) |\Delta c_\varepsilon|^{2}dx+4 \int_{\Omega}(n_\varepsilon+\varepsilon)^{-1}|\nabla n_\varepsilon|^2|\nabla c_\varepsilon|^2 dx
\\
&\le -\frac34\int_{\Omega}n_\varepsilon |\Delta c_\varepsilon|^{2}dx+\frac{\varepsilon}4\int_{\Omega}|\Delta c_\varepsilon|^{2}dx
+\frac12\int_{\Omega}|\nabla\Delta c_\varepsilon|^{2} dx+ \int_{\Omega}(|\nabla u_\varepsilon \nabla c_\varepsilon|^2
+|u_\varepsilon  D^2c_\varepsilon |^2) dx
\\
&+\eta\int_{\Omega}\left((n_\varepsilon+\varepsilon)^{m-2} |\Delta n_\varepsilon|^{2}+(n_\varepsilon+\varepsilon)^{m-4} |\nabla n_\varepsilon|^{4}\right)dx+C_\eta \int_{\Omega}\left((n_\varepsilon+\varepsilon)^{2-m} |\Delta c_\varepsilon|^{2}+(n_\varepsilon+\varepsilon)^{2-m} |\nabla c_\varepsilon|^{4} \right)dx
\\
&\le -\frac34\int_{\Omega}n_\varepsilon |\Delta c_\varepsilon|^{2}dx+\frac{\varepsilon}4\int_{\Omega}|\Delta c_\varepsilon|^{2}dx
+\frac12\int_{\Omega}|\nabla\Delta c_\varepsilon|^{2} dx+ \int_{\Omega}(|\nabla u_\varepsilon \nabla c_\varepsilon|^2
+|u_\varepsilon  D^2c_\varepsilon |^2) dx
\\
&+\eta\int_{\Omega}\left((n_\varepsilon+\varepsilon)^{m-2} |\Delta n_\varepsilon|^{2}+(n_\varepsilon+\varepsilon)^{m-4} |\nabla n_\varepsilon|^{4}\right)dx+C_\eta \int_{\Omega}\left((n_\varepsilon+\varepsilon)^{2-m} |\Delta c_\varepsilon|^{2}+(n_\varepsilon+\varepsilon)^{2-m} |\nabla c_\varepsilon|^{4} \right)dx,
\end{align*}
that is
\allowdisplaybreaks
\begin{align}
&\frac{d}{dt}\int_\Omega |\Delta c_\varepsilon|^{2} dx
+\int_{\Omega}|\nabla\Delta c_\varepsilon|^{2} dx+\frac32\int_{\Omega}n_\varepsilon |\Delta c_\varepsilon|^{2}dx\nonumber
\\
&\le \frac{\varepsilon}2\int_{\Omega}|\Delta c_\varepsilon|^{2}dx
+2 \int_{\Omega}(|\nabla u_\varepsilon \nabla c_\varepsilon|^2
+|u_\varepsilon  D^2c_\varepsilon |^2) dx+2\eta\int_{\Omega}\left((n_\varepsilon+\varepsilon)^{m-2} |\Delta n_\varepsilon|^{2}+(n_\varepsilon+\varepsilon)^{m-4} |\nabla n_\varepsilon|^{4}\right)dx\nonumber
\\
\label{3-39}
&+2C_\eta \int_{\Omega}\left((n_\varepsilon+\varepsilon)^{2-m} |\Delta c_\varepsilon|^{2}+(n_\varepsilon+\varepsilon)^{2-m} |\nabla c_\varepsilon|^{4} \right)dx.
\end{align}
Adding \eqref{3-34} to \eqref{3-39} gives
\begin{align*}
&\frac{d}{dt}\int_\Omega \left(|\Delta c_\varepsilon|^{2}+|\nabla c_\varepsilon|^{4}\right) dx
+\int_{\Omega}|\nabla\Delta c_\varepsilon|^{2} dx+\frac32\int_{\Omega}n_\varepsilon |\Delta c_\varepsilon|^{2}dx
+2\int_{\Omega} n_\varepsilon|\nabla c_\varepsilon|^{4}dx
\\
&+\int_{\Omega}|\nabla c_\varepsilon|^{2}| D^2 c_\varepsilon|^2 dx+2
\int_{\Omega}|\nabla c_\varepsilon|^{2}(\nabla|\nabla c_\varepsilon|)^2 dx
\\
&\le \frac{\varepsilon}2\int_{\Omega}|\Delta c_\varepsilon|^{2}dx
+2 \int_{\Omega}(|\nabla u_\varepsilon \nabla c_\varepsilon|^2
+|u_\varepsilon  D^2c_\varepsilon |^2) dx+2\eta\int_{\Omega}(n_\varepsilon+\varepsilon)^{m-2} |\Delta n_\varepsilon|^{2}dx
\\
&+2C_\eta \int_{\Omega}\left((n_\varepsilon+\varepsilon)^{2-m} |\Delta c_\varepsilon|^{2}+(n_\varepsilon+\varepsilon)^{2-m} |\nabla c_\varepsilon|^{4} \right)dx+(2\eta+\tilde\eta)\int_{\Omega}(n_\varepsilon+\varepsilon)^{m-4}|\nabla n_\varepsilon|^{4}dx
\\
&+
C_{\tilde\eta}\int_{\Omega}|\nabla c_\varepsilon|^{4}dx+C\int_{\Omega}|u_\varepsilon|^6 dx
+C
\\
&\le
2 \int_{\Omega}(|\nabla u_\varepsilon \nabla c_\varepsilon|^2
+|u_\varepsilon  D^2c_\varepsilon |^2) dx+2\eta\int_{\Omega}(n_\varepsilon+\varepsilon)^{m-2} |\Delta n_\varepsilon|^{2}dx+(2\eta+\tilde\eta)\int_{\Omega}(n_\varepsilon+\varepsilon)^{m-4}|\nabla n_\varepsilon|^{4}dx
\\
&+\int_{\Omega}(\frac12 n_\varepsilon |\Delta c_\varepsilon|^{2}+n_\varepsilon|\nabla c_\varepsilon|^{4})dx
+\hat C_\eta \int_{\Omega}\left(|\Delta c_\varepsilon|^{2}+ |\nabla c_\varepsilon|^{4} \right)dx
+C_{\tilde\eta}\int_{\Omega}|\nabla c_\varepsilon|^{4}dx+C\int_{\Omega}|u_\varepsilon|^6 dx+C
\end{align*}
since $0<2-m<1$. By the arbitrariness of $\eta$ and $\tilde\eta$. The above inequality is equivalent to
\begin{align}
&\frac{d}{dt}\int_\Omega \left(|\Delta c_\varepsilon|^{2}+|\nabla c_\varepsilon|^{4}\right) dx
+\int_{\Omega}|\nabla\Delta c_\varepsilon|^{2} dx+\int_{\Omega}n_\varepsilon |\Delta c_\varepsilon|^{2}dx
+\int_{\Omega} n_\varepsilon|\nabla c_\varepsilon|^{4}dx\nonumber
\\
&+\int_{\Omega}|\nabla c_\varepsilon|^{2}| D^2 c_\varepsilon|^2 dx+2
\int_{\Omega}|\nabla c_\varepsilon|^{2}(\nabla|\nabla c_\varepsilon|)^2 dx \nonumber
\\
&\le
2 \int_{\Omega}(|\nabla u_\varepsilon \nabla c_\varepsilon|^2
+|u_\varepsilon  D^2c_\varepsilon |^2) dx+\eta\int_{\Omega}\left((n_\varepsilon+\varepsilon)^{m-2} |\Delta n_\varepsilon|^{2}+(n_\varepsilon+\varepsilon)^{m-4}|\nabla n_\varepsilon|^{4}\right)dx \nonumber
\\
\label{3-40}
&+C_\eta \int_{\Omega}\left(|\Delta c_\varepsilon|^{2}+ |\nabla c_\varepsilon|^{4} \right)dx
+C\int_{\Omega}|u_\varepsilon|^6 dx+C
\end{align}
for any small $\eta>0$.
By Lemma \ref{lem2-1},
$$
\|\nabla c_\varepsilon\|_{W^{2,p}}\le C\|\Delta c_\varepsilon\|_{W^{1,p}}.
$$
Then by Gagliardo-Nirenberg interpolation inequality and \eqref{3-3}, we see that
$$
\| D^2 c_\varepsilon\|_{L^3}^3\le
C_1\|\nabla^3 c_\varepsilon\|_{L^2}^{2}\|c_\varepsilon\|_{L^\infty}+C_2\|c_\varepsilon\|_{L^\infty}^3
\le C_3\left(1+\|\nabla\Delta c_\varepsilon\|_{L^2}^{2}+\|\Delta c_\varepsilon\|_{L^2}^{2}\right).
$$
Thus we have
 \begin{align}
&2 \int_{\Omega}(|\nabla u_\varepsilon \nabla c_\varepsilon|^2
+|u_\varepsilon  D^2c_\varepsilon |^2) dx\nonumber
\\
&\le \eta_1\int_{\Omega}\left(|\nabla c_\varepsilon|^6+| D^2c_\varepsilon |^3\right) dx
+C_{\eta_1}\int_{\Omega}\left(|\nabla u_\varepsilon|^3+|u_\varepsilon|^6\right)dx\nonumber
\\
\label{3-41}
&\le \eta_1\int_{\Omega}\left(|\nabla c_\varepsilon|^6+C_3|\nabla\Delta c_\varepsilon |^2\right) dx
+C_{\eta_1}\int_{\Omega}\left(|\nabla u_\varepsilon|^3+|u_\varepsilon|^6\right)dx
+C_3\eta_1\left(1+\|\Delta c_\varepsilon\|_{L^2}^{2}\right).
 \end{align}
Recalling \eqref{3-37}, and substituting \eqref{3-41} into \eqref{3-40} with $\eta_1$ appropriately small gives \eqref{3-38},
and we complete the proof of this lemma. \hfill $\Box$

\begin{lemma}
\label{lem3-9}
Assume that $1<m\le\frac54$.
Let $(n_\varepsilon, c_\varepsilon, u_\varepsilon, \pi_\varepsilon)$ be the classical
solution of \eqref{3-1} in $(0, T_{\max})$. Then  there exists a constant $C$ such that
\begin{align}
&\sup_{0<t<T_{\max}}\int_\Omega \left(\frac{|\nabla n_\varepsilon|^2}{n_\varepsilon+\varepsilon}+|\Delta c_\varepsilon|^{2}+|\nabla c_\varepsilon|^{4}\right)dx
+\sup_{\tau<t<T_{\max}}\int_{t-\tau}^t\int_\Omega (n_\varepsilon+\varepsilon)| D^2 (n_\varepsilon+\varepsilon)^{\frac{m-1}2}|^2dxds  \nonumber
\\
&+\sup_{\tau<t<T_{\max}}\int_{t-\tau}^t\int_\Omega\left((n_\varepsilon+\varepsilon)^{m-4}|\nabla n_\varepsilon|^4+ (n_\varepsilon+\varepsilon)^{m-2}|\Delta n_\varepsilon|^2\right)dxds  \nonumber
\\
&+\sup_{\tau<t<T_{\max}}\int_{t-\tau}^t\int_\Omega\left(|\nabla\Delta c_\varepsilon|^{2}+n_\varepsilon|\nabla c_\varepsilon|^{4}+n_\varepsilon |\Delta c_\varepsilon|^{2}+|\nabla c_\varepsilon|^{2}\left(| D^2 c_\varepsilon|^2+(\nabla|\nabla c_\varepsilon|)^2\right)\right)dxds
\nonumber
\\
\label{3-42}
&\le C,
\end{align}
where  $C$ is independent of $\varepsilon$ and $T_{\max}$.
\end{lemma}

{\bfseries Proof.} Using Gagliardo-Nirenberg interpolation inequality, we get that
$$
\|n_\varepsilon+\varepsilon\|_{L^{\frac{7m}{3}}}^{\frac{7m}{3}}=\|(n_\varepsilon+\varepsilon)^{\frac m4}\|_{L^{\frac{28}{3}}}^{\frac{28}{3}}
\le C\|(n_\varepsilon+\varepsilon)^{\frac m4}\|_{L^4}^{\frac{16}3}\|\nabla(n_\varepsilon+\varepsilon)^{\frac m4}\|_{L^4}^{4}+C\|(n_\varepsilon+\varepsilon)^{\frac m4}\|_{L^4}^{\frac{28}3}.
$$
Then from \eqref{3-17}, we infer  that
\begin{align}\label{3-43}
\|n_\varepsilon+\varepsilon\|_{L^{\frac{7m}{3}}}^{\frac{7m}{3}}\le C\left(1+\int_\Omega(n_\varepsilon+\varepsilon)^{m-4}|\nabla n_\varepsilon|^{4} dx\right ).
\end{align}
For $u_\varepsilon$,  by \eqref{3-17}, we obtain that
\begin{align}
\sup_{t\in(0,T_{\max})}\|u_\varepsilon(\cdot, t)\|_{L^3}\le &\sup_{t\in(0,T_{\max})}e^{-t}\|u_{\varepsilon 0}\|_{L^3}
+\sup_{t\in(0,T_{\max})}\int_0^t \|e^{-(t-s)A}P(n_\varepsilon(s)\nabla\varphi_\varepsilon(s))\|_{L^3}ds\nonumber
\\
\le &\sup_{t\in(0,T_{\max})}e^{-t}\|u_{\varepsilon 0}\|_{L^3}+
\sup_{t\in(0,T_{\max})}\int_0^t e^{-\lambda(t-s)}(t-s)^{-\frac32(\frac1m-\frac13)}\|n_\varepsilon(s)\nabla\varphi_\varepsilon\|_{L^m}ds
\nonumber
\\
\le &\sup_{t\in(0,T_{\max})}e^{-t}\|u_{\varepsilon 0}\|_{L^3}+\sup_{t\in(0,T_{\max})}
\|n_\varepsilon\|_{L^m}\|\nabla\varphi_\varepsilon\|_{L^\infty}
\int_0^{\infty} e^{-\lambda s}s^{-\frac32(\frac1m-\frac13)}ds\nonumber
\\
\label{3-44}
\le & C
\end{align}
since $0<\frac32(\frac1m-\frac13)<1$.
Using Gagliardo-Nirenberg  inequality and \eqref{3-44}, we further obtain that
\begin{align}\label{3-45}
\|u_\varepsilon\|_{L^{7m}}^{7m}\le C\|u_\varepsilon\|_{L^3}^{\frac{14m}{3}}\|A u_\varepsilon\|_{L^{\frac{7m}3}}^{\frac{7m}3}
+C \|u_\varepsilon\|_{L^3}^{7m}\le C\left(1+\|A u_\varepsilon\|_{L^{\frac{7m}3}}^{\frac{7m}3}\right),
\\
\label{3-46}
\|\nabla u_\varepsilon\|_{L^{\frac{7m}2}}^{\frac{7m}2}\le  C\|u_\varepsilon\|_{L^3}^{\frac{7m}{6}}\|A u_\varepsilon\|_{L^{\frac{7m}3}}^{\frac{7m}3}
+C \|u_\varepsilon\|_{L^3}^{\frac{7m}2}\le C\left(1+\|A u_\varepsilon\|_{L^{\frac{7m}3}}^{\frac{7m}3}\right).
\end{align}
By  $L^{p, q}$ theory of Stokes operator \cite{FK, HP} and \eqref{3-43}, we also have
\begin{align}
\sup_{t\in(\tau,T_{\max})}\int_{t-\tau}^{t}
\|A u_\varepsilon\|_{L^{\frac{7m}3}}^{\frac{7m}3}ds&\le C_1+C_2\sup_{t\in(\tau,T_{\max})}\int_{t-\tau}^{t}
\|n_\varepsilon\|_{L^{\frac{7m}3}}^{\frac{7m}3}ds\nonumber
\\
\label{3-47}
&\le  C_3+C_4\sup_{t\in(\tau,T_{\max})}\int_{t-\tau}^{t}\int_\Omega(n_\varepsilon+\varepsilon)^{m-4}|\nabla n_\varepsilon|^{4} dxds.
\end{align}
Combining \eqref{3-45}, \eqref{3-46} and \eqref{3-47}, we derive that
\begin{align}\label{3-48}
\sup_{t\in(\tau,T_{\max})}\int_{t-\tau}^{t}\int_\Omega
\left(|u_\varepsilon|^{7m}+|\nabla u_\varepsilon|^{\frac{7m}2}\right)dxds  \le  C+C\sup_{t\in(\tau,T_{\max})}\int_{t-\tau}^{t}\int_\Omega(n_\varepsilon+\varepsilon)^{m-4}|\nabla n_\varepsilon|^{4} dxds.
\end{align}
Adding $\int_\Omega\frac{|\nabla n_\varepsilon|^2}{n_\varepsilon+\varepsilon} dx$ on both sides of \eqref{3-20}, combining the resulting inequality
and  \eqref{3-42}, and noting that $(n_\varepsilon+\varepsilon)^{-1}\le
(n_\varepsilon+\varepsilon)^{2m-4}+(n_\varepsilon+\varepsilon)^{m-2}$ since $2m-4<-1<m-2$, then we see that
\begin{align}
&\frac{d}{dt}\int_\Omega \left(\frac{|\nabla n_\varepsilon|^2}{n_\varepsilon+\varepsilon}+|\Delta c_\varepsilon|^{2}+|\nabla c_\varepsilon|^{4}\right)dx
+\frac{2m}{(m-1)^2}\int_\Omega (n_\varepsilon+\varepsilon)| D^2 (n_\varepsilon+\varepsilon)^{\frac{m-1}2}|^2 dx+\frac12\int_{\Omega}|\nabla\Delta c_\varepsilon|^{2}dx\nonumber
\\
&
+\int_{\Omega} \left(n_\varepsilon|\nabla c_\varepsilon|^{4}+n_\varepsilon |\Delta c_\varepsilon|^{2}\right)dx+\int_{\Omega}\frac12|\nabla c_\varepsilon|^{2}\left(| D^2 c_\varepsilon|^2+(\nabla|\nabla c_\varepsilon|)^2\right) dx+\int_\Omega \left(|\Delta c_\varepsilon|^{2}+|\nabla c_\varepsilon|^{4}+\frac{|\nabla n_\varepsilon|^2}{n_\varepsilon+\varepsilon}\right) dx\nonumber
\\
&\le\int_\Omega\frac{|\nabla n_\varepsilon|^2}{n_\varepsilon+\varepsilon}dx+ C \int_{\Omega}(n_\varepsilon+\varepsilon)^{m-2}|\nabla n_\varepsilon |^2 dx+C\int_\Omega(n_\varepsilon+\varepsilon)^{4-3m} |u_{\varepsilon}|^4dx\nonumber
\\
&+C\int_\Omega\left((n_\varepsilon+\varepsilon)^{4-3m} |\nabla c_{\varepsilon}|^4+(n_\varepsilon+\varepsilon)^{2-m}|\Delta c_\varepsilon|^2\right) dx+
\eta\int_{\Omega}\left((n_\varepsilon+\varepsilon)^{m-2} |\Delta n_\varepsilon|^{2}+(n_\varepsilon+\varepsilon)^{m-4}|\nabla n_\varepsilon|^{4}\right)dx
\nonumber
\\
&+C_\eta \int_{\Omega}\left(|\Delta c_\varepsilon|^{2}+ |\nabla c_\varepsilon|^{4} \right)dx+C\int_{\Omega}\left(|\nabla u_\varepsilon|^3+|u_\varepsilon|^6\right)dx+C\nonumber
\\
&\le\int_\Omega (n_\varepsilon+\varepsilon)^{2m-4}|\nabla n_\varepsilon|^2dx+ \tilde C \int_{\Omega}(n_\varepsilon+\varepsilon)^{m-2}|\nabla n_\varepsilon |^2 dx+C\int_\Omega(n_\varepsilon+\varepsilon)^{4-3m} |u_{\varepsilon}|^4dx\nonumber
\\
&+\frac12\int_\Omega\left((n_\varepsilon+\varepsilon)|\nabla c_{\varepsilon}|^4+(n_\varepsilon+\varepsilon)|\Delta c_\varepsilon|^2\right) dx+
\eta\int_{\Omega}\left((n_\varepsilon+\varepsilon)^{m-2} |\Delta n_\varepsilon|^{2}+(n_\varepsilon+\varepsilon)^{m-4}|\nabla n_\varepsilon|^{4}\right)dx
\nonumber
\\
\label{3-49}
&+\tilde C_\eta \int_{\Omega}\left(|\Delta c_\varepsilon|^{2}+ |\nabla c_\varepsilon|^{4} \right)dx+C\int_{\Omega}\left(|\nabla u_\varepsilon|^3+|u_\varepsilon|^6\right)dx+C
\end{align}
for any small $\eta>0$,
since $0<4-3m<1$, $0<2-m<1$. Recalling \eqref{3-27} and \eqref{3-30}, taking $\eta$ appropriately small in \eqref{3-49}, we obtain that
\begin{align}
&\frac{d}{dt}\int_\Omega \left(\frac{|\nabla n_\varepsilon|^2}{n_\varepsilon+\varepsilon}+|\Delta c_\varepsilon|^{2}+|\nabla c_\varepsilon|^{4}\right)dx
+\frac{3m}{2(m-1)^2}\int_\Omega (n_\varepsilon+\varepsilon)| D^2 (n_\varepsilon+\varepsilon)^{\frac{m-1}2}|^2 dx+\frac12\int_{\Omega}|\nabla\Delta c_\varepsilon|^{2}dx\nonumber
\\
&
+\frac12\int_{\Omega} \left(n_\varepsilon|\nabla c_\varepsilon|^{4}+n_\varepsilon |\Delta c_\varepsilon|^{2}\right)dx+\frac12\int_{\Omega}|\nabla c_\varepsilon|^{2}\left(| D^2 c_\varepsilon|^2+(\nabla|\nabla c_\varepsilon|)^2\right) dx\nonumber
\\
&+\int_\Omega \left(|\Delta c_\varepsilon|^{2}+|\nabla c_\varepsilon|^{4}+\frac{|\nabla n_\varepsilon|^2}{n_\varepsilon+\varepsilon}\right) dx\nonumber
\\
&\le\int_\Omega (n_\varepsilon+\varepsilon)^{2m-4}|\nabla n_\varepsilon|^2dx+ \tilde C \int_{\Omega}(n_\varepsilon+\varepsilon)^{m-2}|\nabla n_\varepsilon |^2 dx+C\int_\Omega(n_\varepsilon+\varepsilon)^{4-3m} |u_{\varepsilon}|^4dx\nonumber
\\
\label{3-50}
&+C \int_{\Omega}\left(|\Delta c_\varepsilon|^{2}+ |\nabla c_\varepsilon|^{4} \right)dx+C\int_{\Omega}\left(|\nabla u_\varepsilon|^3+|u_\varepsilon|^6\right)dx+C.
\end{align}
Recalling \eqref{3-43} and \eqref{3-26}, for any small constant $\sigma>0$, we get that
\begin{align}
&C\int_\Omega(n_\varepsilon+\varepsilon)^{4-3m} |u_{\varepsilon}|^4dx\le \sigma\int_\Omega(n_\varepsilon+\varepsilon)^{\frac{7m}3}dx+
C_{\sigma}\int_\Omega |u_{\varepsilon}|^{\frac{7m}{4m-3}}dx
\nonumber
\\
&\le \sigma \tilde C\left(1+\int_\Omega(n_\varepsilon+\varepsilon)^{m-4}|\nabla n_\varepsilon|^{4} dx\right )
+C_{\sigma}\int_\Omega |u_{\varepsilon}|^{\frac{7m}{4m-3}}dx\nonumber
\\
\label{3-51}
&\le \frac{m}{2(m-1)^2}\int_\Omega (n_\varepsilon+\varepsilon)| D^2 (n_\varepsilon+\varepsilon)^{\frac{m-1}2}|^2 dx
+\hat C\int_\Omega |u_{\varepsilon}|^{\frac{7m}{4m-3}}dx+\hat C.
 \end{align}
Substituting it into \eqref{3-50} yields
\begin{align}
&\frac{d}{dt}\int_\Omega \left(\frac{|\nabla n_\varepsilon|^2}{n_\varepsilon+\varepsilon}+|\Delta c_\varepsilon|^{2}+|\nabla c_\varepsilon|^{4}\right)dx
+\frac{m}{(m-1)^2}\int_\Omega (n_\varepsilon+\varepsilon)| D^2 (n_\varepsilon+\varepsilon)^{\frac{m-1}2}|^2 dx+\frac12\int_{\Omega}|\nabla\Delta c_\varepsilon|^{2}dx\nonumber
\\
&
+\frac12\int_{\Omega} \left(n_\varepsilon|\nabla c_\varepsilon|^{4}+n_\varepsilon |\Delta c_\varepsilon|^{2}\right)dx+\frac12\int_{\Omega}|\nabla c_\varepsilon|^{2}\left(| D^2 c_\varepsilon|^2+(\nabla|\nabla c_\varepsilon|)^2\right) dx\nonumber
\\
&+\int_\Omega \left(|\Delta c_\varepsilon|^{2}+|\nabla c_\varepsilon|^{4}+\frac{|\nabla n_\varepsilon|^2}{n_\varepsilon+\varepsilon}\right) dx\nonumber
\\
&\le\int_\Omega (n_\varepsilon+\varepsilon)^{2m-4}|\nabla n_\varepsilon|^2dx+ \tilde C \int_{\Omega}(n_\varepsilon+\varepsilon)^{m-2}|\nabla n_\varepsilon |^2 dx+C \int_{\Omega}\left(|\Delta c_\varepsilon|^{2}+ |\nabla c_\varepsilon|^{4} \right)dx\nonumber
\\
\label{3-52}
&+C\int_{\Omega}\left(|\nabla u_\varepsilon|^3+|u_\varepsilon|^6+|u_{\varepsilon}|^{\frac{7m}{4m-3}}\right)dx+C.
\end{align}
Using \eqref{3-5}, \eqref{3-33}, \eqref{3-48} and \eqref{3-26}, we further have
\begin{align}
&\sup_{0<t<T_{\max}}\int_\Omega \left(\frac{|\nabla n_\varepsilon|^2}{n_\varepsilon+\varepsilon}+|\Delta c_\varepsilon|^{2}+|\nabla c_\varepsilon|^{4}\right)dx
+\frac{m}{(m-1)^2}\sup_{\tau<t<T_{\max}}\int_{t-\tau}^t\int_\Omega (n_\varepsilon+\varepsilon)| D^2 (n_\varepsilon+\varepsilon)^{\frac{m-1}2}|^2dxds  \nonumber
\\
&+\sup_{\tau<t<T_{\max}}\int_{t-\tau}^t\int_\Omega\left(|\nabla\Delta c_\varepsilon|^{2}+n_\varepsilon|\nabla c_\varepsilon|^{4}+n_\varepsilon |\Delta c_\varepsilon|^{2}+|\nabla c_\varepsilon|^{2}\left(| D^2 c_\varepsilon|^2+(\nabla|\nabla c_\varepsilon|)^2\right)\right)dxds
\nonumber
\\
&\le C\sup_{\tau<t<T_{\max}}\int_{t-\tau}^t\int_\Omega\left( (n_\varepsilon+\varepsilon)^{2m-4}|\nabla n_\varepsilon|^2+|\nabla (n_\varepsilon+\varepsilon)^{m/2} |^2+\Delta c_\varepsilon|^{2}+ |\nabla c_\varepsilon|^{4}\right)dxds\nonumber
\\
&+C\sup_{\tau<t<T_{\max}}\int_{t-\tau}^t\int_{\Omega}\left(|\nabla u_\varepsilon|^3+|u_\varepsilon|^6+|u_{\varepsilon}|^{\frac{7m}{4m-3}}\right)dxds+C\nonumber
\\
&\le C_\eta+\eta\sup_{\tau<t<T_{\max}}\int_{t-\tau}^t\int_{\Omega}\left(|\nabla u_\varepsilon|^{\frac{7m}2}+|u_\varepsilon|^{7m}\right)dxds\nonumber
\\
\label{3-53}
&\le \hat C+\frac{m}{2(m-1)^2}\sup_{\tau<t<T_{\max}}\int_{t-\tau}^t\int_\Omega (n_\varepsilon+\varepsilon)| D^2 (n_\varepsilon+\varepsilon)^{\frac{m-1}2}|^2dxds,
\end{align}
since $\max\{6, \frac{7m}{4m-3}\}<7m$, for any small $\eta>0$.
Combining with \eqref{3-26} and \eqref{3-30}, we complete the proof.
\hfill $\Box$

\begin{lemma}
\label{lem3-10}
Assume that $1<m\le\frac54$.
Let $(n_\varepsilon, c_\varepsilon, u_\varepsilon, \pi_\varepsilon)$ be the classical
solution of \eqref{3-1} in $(0, T_{\max})$. Then  for any $p\ge 2$,
\begin{align}
\label{3-54}
&\sup_{0<t<T_{\max}}\left(\|u_\varepsilon(\cdot, t)\|_{W^{1,\infty}}+\|A^\beta u_\varepsilon(\cdot, t)\|_{L^{p}}\right) \le C_,
\\
\label{3-55}
&\sup_{0<t<T_{\max}}\left(\|c_\varepsilon(\cdot, t)\|_{W^{1,\infty}}+\|n_\varepsilon(\cdot, t)\|_{L^{\infty}}\right)\le C,
\end{align}
where  $C_1$, $C_2$ are independent of $\varepsilon$ and $T_{\max}$.
\end{lemma}

{\bfseries Proof.}
By Sobolev imbedding inequalities, and using  Lemma \ref{lem3-9}, we have
$$
\|n_\varepsilon+\varepsilon\|_{L^3}=\|\sqrt{n_\varepsilon+\varepsilon}\|_{L^6}\le C \|\sqrt{n_\varepsilon+\varepsilon}\|_{H^1}\le C_1,
$$
$$
\|\nabla c_\varepsilon\|_{L^6}\le C(\|\nabla c_\varepsilon\|_{L^4}+\|\Delta c_\varepsilon\|_{L^2})\le C_2.
$$
Recalling \eqref{3-19}, for any $r>m$,
\begin{align*}
&\frac1{r+1}\frac{d}{dt}\int_\Omega |n_\varepsilon+\varepsilon|^{r+1} dx
+\frac{rm}2\int_\Omega |n_\varepsilon+\varepsilon|^{m+r-2}|\nabla n_\varepsilon|^2 dx+\int_\Omega |n_\varepsilon+\varepsilon|^{r+1} dx
\\
\le & C\int_\Omega |n_\varepsilon+\varepsilon|^{r+2-m}|\nabla c_\varepsilon|^{2}dx+C
\\
\le &C\|\nabla c_\varepsilon\|_{L^6}^2\|n_\varepsilon+\varepsilon\|_{L^{\frac{3(r+2-m)}2}}^{r+2-m}+C
\\
\le &C_3\|(n_\varepsilon+\varepsilon)^{\frac{r+m}2}\|_{L^{\frac{3(r+2-m)}{r+m}}}^{\frac{2(r+2-m)}{r+m}}+C
\\
\le &C_3\|(n_\varepsilon+\varepsilon)^{\frac{r+m}2}\|_{L^{\frac 6{r+m}}}^{\frac{2(r+2-m)}{r+m}-\frac{2(r-m)}{r+m-1}}
\|\nabla(n_\varepsilon+\varepsilon)^{\frac{r+m}2}\|_{L^2}^{\frac{2(r-m)}{r+m-1}}+C_4
\\
\le & \frac{rm}4\int_\Omega |n_\varepsilon+\varepsilon|^{m+r-2}|\nabla n_\varepsilon|^2 dx+C_5,
\end{align*}
which implies that for any $r>m$.
\begin{align}\label{3-57}
\sup_{0<t<T_{\max}}\int_\Omega |n_\varepsilon+\varepsilon|^{r+1} dx+\sup_{\tau<t<T_{\max}}\int_{t-\tau}^t\int_\Omega|n_\varepsilon+\varepsilon|^{m+r-2}|\nabla n_\varepsilon|^2 dx ds\le C.
\end{align}
For any $\beta\in(\frac34, 1)$, and any $p\ge 2$,
\begin{align*}
\|A^{\beta}u_\varepsilon\|_{L^p}\le &e^{-t}\|A^{\beta}u_{\varepsilon 0}\|_{L^p}
+\int_0^t \|A^{\beta}e^{-(t-s)A}P(n_\varepsilon(s)\nabla\varphi_\varepsilon(s))\|_{L^p}ds
\\
\le &e^{-t}\|A^{\beta}u_{\varepsilon 0}\|_{L^p}+\int_0^t e^{-\lambda(t-s)}(t-s)^{-\beta}\|n_\varepsilon(s)
\nabla\varphi_\varepsilon\|_{L^p}ds
\\
\le &e^{-t}\|A^{\beta}u_{\varepsilon 0}\|_{L^p}+\int_0^t e^{-\lambda(t-s)}(t-s)^{-\beta}
\|n_\varepsilon\|_{L^p}\|\nabla\varphi_\varepsilon\|_{L^\infty}ds
\\
\le &C,
\end{align*}
by embedding theorem, $u_\varepsilon, \nabla u_\varepsilon\in L^\infty(\Omega\times (0, T_{\max}))$ since $\beta>\frac34$, $\forall p\ge 2$,
and \eqref{3-54} is proved.

For $c_\varepsilon$,  we have
\begin{align*}
\|\nabla c_\varepsilon\|_{L^\infty}\le & e^{-t}\|\nabla c_{\varepsilon 0}\|_{L^\infty}+\int_0^t e^{-(t-s)}
\left\|\nabla\left(e^{(t-s)\Delta}\Big(-u_\varepsilon\cdot\nabla n_\varepsilon-c_\varepsilon n_\varepsilon)\Big)\right)\right\|_{L^\infty}ds
\\
\le & e^{-t}\|\nabla c_{\varepsilon 0}\|_{L^\infty}+\int_0^t e^{-(t-s)}(t-s)^{-\frac12-\frac14}
\left\|-u_\varepsilon\cdot\nabla c_\varepsilon-c_\varepsilon n_\varepsilon\right\|_{L^6}ds
\\
\le & e^{-t}\|\nabla c_{\varepsilon 0}\|_{L^\infty}+\int_0^t e^{-(t-s)}(t-s)^{-\frac34}
(\|u_\varepsilon \|_{L^\infty}\|\nabla c_\varepsilon\|_{L^6}+\|c_\varepsilon\|_{L^\infty}\|n_\varepsilon\|_{L^6})ds
\\
\le & e^{-t}\|\nabla c_{\varepsilon 0}\|_{L^\infty}+C\int_0^t e^{-(t-s)}(t-s)^{-\frac34}ds
\\
\le & e^{-t}\|\nabla c_{\varepsilon 0}\|_{L^\infty}+C\int_0^t e^{-s}s^{-\frac34}ds
\\
\le & \tilde C, \qquad  \text{for any $0<t<T_{\max}$}.
\end{align*}
Next, by a standard Moser iteration technique, we obtain the $L^\infty$-norm estimate of $n_\varepsilon$, and we complete the proof of this lemma.
 \hfill $\Box$

\begin{remark}
\label{re-1}
Although we assume that $1<m\le \frac54$ in Lemma \ref{lem3-10},  in fact, for any $m>1$,  the following inequality holds
\begin{align*}
&\sup_{0<t<T_{\max}}\left(\|c_\varepsilon(\cdot, t)\|_{W^{1,\infty}}+\|u_\varepsilon(\cdot, t)\|_{W^{1,\infty}}+\|n_\varepsilon(\cdot, t)\|_{L^{\infty}}\right)\le C.
\end{align*}
\end{remark}

\medskip

 Then using Lemma \ref{lem3-1}, $T_{\max}=+\infty$, that is, the classical solution of \eqref{3-1} exists globally for any $m>1$.
 Next, we give some higher regularity estimates.

\begin{lemma}
\label{lem3-11}
Assume that $m>1$.
Let $(n_\varepsilon, c_\varepsilon, u_\varepsilon, \pi_\varepsilon)$ be the global classical
solution of \eqref{3-1}. Then
\begin{align}\label{3-58}
&\sup_{t>0}\int_t^{t+1}\left(\|u_{\varepsilon t}\|_{L^p}^p+\|u_\varepsilon\|_{W^{2,p}}^p+\|\nabla\pi_\varepsilon\|_{L^p}^p\right)ds\le C, \quad \text{for any $p>1$},
\\
\label{3-59}
&\sup_{t>0}\int_t^{t+1}\left(\|c_{\varepsilon t}\|_{L^p}^p+\|c_\varepsilon\|_{W^{2,p}}^p\right)ds\le C, \quad \text{for any $p>1$},
\\
\label{3-60}
&\sup_{t\in(0, +\infty)}\int_\Omega |\nabla (n_\varepsilon+\varepsilon)^m|^2dx
+\sup_{t\in(0, +\infty)}\int_t^{t+1}\int_\Omega (n_\varepsilon+\varepsilon)^{m-1} \left|\frac{\partial n_\varepsilon}{\partial t}\right|^2 dx\le C.
\end{align}
In particular, if $1<m\le \frac54$, we also have
\begin{align}
\label{3-61}
&\sup_{t>0}\int_t^{t+1}\int_\Omega \left(\left|\frac{\partial n_\varepsilon}{\partial t}\right|^2+|\Delta (n_\varepsilon+\varepsilon)^m|^2\right) dxds\le C.
\end{align}
Here, all these constants $C$ are independent of $\varepsilon$.
\end{lemma}

{\bfseries Proof.} By Lemma \ref{lem2-3} and the $L^{p,q}$ theory of Stokes operator \cite{FK, HP},
\eqref{3-58} and \eqref{3-59} are proved.

Multiplying the first equation of \eqref{3-1} by
$\frac{\partial(n_\varepsilon+\varepsilon)^m)}{\partial t}$,
and integrating it over $\Omega$ gives
\allowdisplaybreaks
\begin{align*}
&\frac1{2}\frac{d}{dt}\int_\Omega |\nabla (n_\varepsilon+\varepsilon)^m|^2dx
+m\int_\Omega (n_\varepsilon+\varepsilon)^{m-1} \left|\frac{\partial n_\varepsilon}{\partial t}\right|^2 dx
+\int_\Omega |\nabla (n_\varepsilon+\varepsilon)^m|^2
\\
&= -m\chi\int_\Omega (n_\varepsilon+\varepsilon)^{m-1} \frac{\partial n_\varepsilon}{\partial t} \nabla\cdot(n_{\varepsilon}\nabla c_{\varepsilon})
dx-m\int_\Omega (n_\varepsilon+\varepsilon)^{m-1} \frac{\partial n_\varepsilon}{\partial t} u_\varepsilon\cdot\nabla n_\varepsilon dx+\int_\Omega |\nabla (n_\varepsilon+\varepsilon)^m|^2
\\
&\le  m\chi^2\int_\Omega|\nabla\cdot(n_{\varepsilon}\nabla c_{\varepsilon})|^2
(n_\varepsilon+\varepsilon)^{m-1} dx+m\int_\Omega |u_\varepsilon\cdot\nabla n_\varepsilon|^2(n_\varepsilon+\varepsilon)^{m-1}dx
\\
&+\int_\Omega m^2(n_\varepsilon+\varepsilon)^{2(m-1)}|\nabla n_\varepsilon|^2
+\frac{m}2\int_\Omega (n_\varepsilon+\varepsilon)^{m-1} \left|\frac{\partial n_\varepsilon}{\partial t}\right|^2 dx.
\end{align*}
By remark \ref{re-1}, we further have
\begin{align*}
&\frac1{2}\frac{d}{dt}\int_\Omega |\nabla (n_\varepsilon+\varepsilon)^m|^2dx
+\frac m2\int_\Omega (n_\varepsilon+\varepsilon)^{m-1} \left|\frac{\partial n_\varepsilon}{\partial t}\right|^2 dx
+\int_\Omega |\nabla (n_\varepsilon+\varepsilon)^m|^2
\\
&\le C\int_\Omega((n_\varepsilon+\varepsilon)^{m-2}|\nabla n_{\varepsilon}|^2+|\Delta c_{\varepsilon}|^2) dx.
\end{align*}
Recalling \eqref{3-5}, then
$$
\sup_{t\in(0, +\infty)}\int_\Omega |\nabla (n_\varepsilon+\varepsilon)^m|^2dx
+\sup_{t\in(0, +\infty)}\int_t^{t+1}\int_\Omega (n_\varepsilon+\varepsilon)^{m-1} \left|\frac{\partial n_\varepsilon}{\partial t}\right|^2 dx\le C.
$$

Multiplying the first equation of \eqref{3-1} by $\frac{\partial n_\varepsilon}{\partial t}$,  integrating it over $\Omega$,
and using  \eqref{3-54} and \eqref{3-55} gives
\begin{align*}
&\int_\Omega\left|\frac{\partial n_\varepsilon}{\partial t}\right|^2dx+\frac{m}2\frac{d}{dt}\int_\Omega\left((n_\varepsilon+\varepsilon)^{m-1}|\nabla n_\varepsilon|^2\right)dx +\int_\Omega\left((n_\varepsilon+\varepsilon)^{m-1}|\nabla n_\varepsilon|^2\right)dx
\\
&=\int_\Omega\left(m(m-1)(n_\varepsilon+\varepsilon)^{m-2}|\nabla n_\varepsilon|^2- u_\varepsilon\nabla n_\varepsilon -\chi \nabla\cdot(n_\varepsilon\nabla c_\varepsilon) \right) n_{\varepsilon t}dx+\int_\Omega\left((n_\varepsilon+\varepsilon)^{m-1}|\nabla n_\varepsilon|^2\right)dx
\nonumber
\\
&\le \frac12\int_\Omega\left|\frac{\partial n_\varepsilon}{\partial t}\right|^2dx +C\int_\Omega \left((n_\varepsilon+\varepsilon)^{m-4}|\nabla n_\varepsilon|^4+|\nabla n_\varepsilon|^2+|\Delta c_\varepsilon|^2\right)dx.
\end{align*}
Using \eqref{3-42}, then we infer that
$$
\sup_{t>0}\int_\Omega\left((n_\varepsilon+\varepsilon)^{m-1}|\nabla n_\varepsilon|^2\right)dx+\sup_{t>0}\int_t^{t+1}\int_\Omega\left|\frac{\partial n_\varepsilon}{\partial t}\right|^2dxds\le C.
$$
Similarly, multiplying the first equation of \eqref{3-1} by $\Delta (n_\varepsilon+\varepsilon)^m$,  integrating it over $\Omega$,
and using  \eqref{3-42}, \eqref{3-54} and \eqref{3-55}, we also infer that
$$
\sup_{t>0}\int_t^{t+1}\int_\Omega |\Delta (n_\varepsilon+\varepsilon)^m|^2 dxds\le C.
$$
This lemma is proved. \hfill $\Box$

\medskip
{\it \bfseries Proof of Theorem \ref{thm-1}.}
Since $(n_\varepsilon, c _\varepsilon, u_\varepsilon, \pi_\varepsilon)$ is the classical solution of
\eqref{3-1}, then we have
\begin{align*}
&-\iint_{Q_T}n_\varepsilon\psi_t dxdt-\int_\Omega n_\varepsilon(x,0)\psi(x,0)dx
+\iint_{Q_T} \nabla (n_\varepsilon+\varepsilon)^m\nabla\psi dxdt
\\
&=\iint_{Q_T}u_\varepsilon n_\varepsilon \nabla \psi dxdt+
\chi\iint_{Q_T}n_\varepsilon\nabla c_\varepsilon\nabla\psi dxdt
\\
&-\iint_{Q_T}c_\varepsilon\phi_tdxdt-\int_\Omega c_\varepsilon(x,0)\phi(x,0)dx
+\iint_{Q_T}\nabla c_\varepsilon\nabla\phi dxdt=\iint_{Q_T}u_\varepsilon c_\varepsilon \nabla \phi dxdt-\iint_{Q_T}
c_{\varepsilon}n_{\varepsilon}\phi dxdt
\\
&-\iint_{Q_T}u_\varepsilon\Phi_tdxdt-\int_\Omega u_\varepsilon(x,0)\Phi(x,0)dx
+\iint_{Q_T}\nabla u_\varepsilon\nabla\Phi dxdt+\iint_{Q_T}\nabla \pi_\varepsilon \Phi dxdt=\iint_{Q_T}n_\varepsilon \nabla \varphi\Phi dxdt
\end{align*}
for any $\psi, \phi, \Phi\in C^\infty(\overline Q_T)$ with  $\phi(x,T)=0$.
Using \eqref{3-5}, remark \ref{re-1}, Lemma \ref{lem3-11}, and Sobolev compact embedding theorem,
for any $T>0$, letting $\varepsilon\to 0$, we have
\begin{align*}
&n_\varepsilon, \varepsilon+n_\varepsilon \to n,    && \text{in $L^p(Q_T)$, for any $p\in (1,+\infty)$},
\\
&c_\varepsilon \to c, \quad u_\varepsilon \to u,  && \text{uniformly},
\\
&n_\varepsilon \stackrel{*}{\rightharpoonup} n, && \text{in $L^\infty(Q_T)$},
\\
&\nabla (n_\varepsilon+\varepsilon)^m \rightharpoonup \nabla n^m, && \text{in  $L^2(Q_T)$},
\\
&\nabla c_\varepsilon \to \nabla c, \nabla u_\varepsilon \to\nabla u,  && \text{in $L^p(Q_T)$, for any $p\in (1,+\infty)$}.
\end{align*}
which means $(n, c, u, \pi)$ is the global bounded weak solution of \eqref{1-1} with \eqref{1-3}--\eqref{1-4} hold. \hfill $\Box$

\medskip

{\it \bfseries Proof of Theorem \ref{thm-2}.} Using Lemma \ref{lem3-9}, Lemma \ref{lem3-10}, Lemma \ref{lem3-11}, and Sobolev imbedding theorem,
it is easy to see that for any $T>0$,
\begin{align*}
&n_\varepsilon, \quad n_\varepsilon+\varepsilon \to n, && \text{in $L^p(Q_T)$, for any $p\in (1,+\infty)$},
\\
&n_\varepsilon \stackrel{*}{\rightharpoonup} n, && \text{in $L^\infty(Q_T)$},
\\
&\nabla n_\varepsilon \to \nabla n,&&  \text{in $L^p(Q_T)$, for any $p\in (1,6)$},
\\
&n_{\varepsilon t} \rightharpoonup n_t, \quad \Delta (n_\varepsilon+\varepsilon)^m \rightharpoonup \Delta n^m, && \text{in $L^2(Q_T)$},
\\
&c_\varepsilon \to c, \quad u_\varepsilon \to u, &&\text{uniformly in $Q_T$},
\\
&\nabla c_\varepsilon \to \nabla c,\quad \nabla u_\varepsilon \to \nabla u, &&  \text{in $L^p(Q_T)$, for any $p>1$},
\\
&c_{\varepsilon t} \rightharpoonup c_t, \quad \Delta c_\varepsilon \rightharpoonup \Delta c,  && \text{in $L^p(Q_T)$, for any $p>1$},
\\
&u_{\varepsilon t} \rightharpoonup u_t, \quad \Delta u_\varepsilon \rightharpoonup \Delta u,  && \text{in $L^p(Q_T)$, for any $p>1$},
\\
&\nabla\pi_\varepsilon \rightharpoonup \nabla\pi,  && \text{in $L^p(Q_T)$, for any $p>1$}.
\end{align*}
Since $(n_\varepsilon, c _\varepsilon, u_\varepsilon, \pi_\varepsilon)$ is the classical solution of
\eqref{3-1}, then for any given $T>0$,
\begin{align*}
&\iint_{Q_T}\left(\frac{\partial n_\varepsilon}{\partial t}+u_\varepsilon\cdot\nabla n_\varepsilon-\Delta(n_\varepsilon+\varepsilon)^m+\chi\nabla\cdot (n_\varepsilon\nabla c_\varepsilon)\right)\phi dxdt, \qquad \text{for any $\phi\in L^2(Q_T)$,}
\\
&\iint_{Q_T}\left(\frac{\partial c_\varepsilon}{\partial t}+u_\varepsilon\cdot\nabla c_\varepsilon-\Delta c_\varepsilon+n_\varepsilon c_\varepsilon\right)\psi dxdt,
\qquad \text{for any $\psi\in L^q(Q_T)$,}
\\
&\iint_{Q_T}\left(\frac{\partial u_\varepsilon}{\partial t}-\Delta u_\varepsilon+\nabla \pi_\varepsilon-n_\varepsilon \nabla \varphi\right)\Phi dxdt,
\qquad \text{for any $\Phi\in L^q(Q_T)$,}
\end{align*}
for any $q>1$.  Letting $\varepsilon\to 0$, then we arrive at
\begin{align*}
&\iint_{Q_T}\left(\frac{\partial n}{\partial t}+u\cdot\nabla n-\Delta n^m+\chi\nabla\cdot (n\nabla c)\right)\phi dxdt, \qquad \text{for any $\phi\in L^2(Q_T)$,}
\\
&\iint_{Q_T}\left(\frac{\partial c}{\partial t}+u\cdot\nabla c-\Delta c+n c\right)\psi dxdt,
\qquad \text{for any $\psi\in L^q(Q_T)$,}
\\
&\iint_{Q_T}\left(\frac{\partial u}{\partial t}-\Delta u+\nabla \pi-n\nabla \varphi\right)\Phi dxdt,
\qquad \text{for any $\Phi\in L^q(Q_T)$,}
\end{align*}
which means $(n, c, u, \pi)$ is the global strong solution of \eqref{1-1} with \eqref{1-5}--\eqref{1-6} hold. \hfill $\Box$

\end{document}